\documentclass[twocolumn]{autart}
\usepackage{listings}
\usepackage{amsmath,amssymb,pstricks,color,dsfont}
\usepackage{savesym}
\savesymbol{AND}
\usepackage{graphicx}
\usepackage{psfrag,graphicx,epsfig}
\usepackage{algorithm,algorithmic,xspace}
\usepackage{subfig}
\usepackage{float}
\usepackage{graphicx}
\usepackage{epstopdf}
\usepackage{url}

\newcommand{\nnum}{\nonumber}

\newcommand{\EQ}{\begin{eqnarray}}
\newcommand{\EN}{\end{eqnarray}}
\newcommand{\EQQ}{\begin{eqnarray*}}
\newcommand{\ENN}{\end{eqnarray*}}

\newcommand{\bremark}{\begin{remark} \begin{rm} }
\newcommand{\eremark}{ \end{rm} \rule{1mm}{2mm}
\end{remark} }
\newcommand{\btheorem}{\begin{theorem} \begin{rm} }
\newcommand{\etheorem}{ \end{rm} \rule{1mm}{2mm}
\end{theorem} }
\newcommand{\blemma}{\begin{lemma} \begin{rm} }
\newcommand{\elemma}{ \end{rm} \rule{1mm}{2mm}
\end{lemma} }
\newcommand{\bcorollary}{\begin{corollary} \begin{rm} }
\newcommand{\ecorollary}{ \end{rm} \rule{1mm}{2mm}
\end{corollary} }
\newcommand{\bdefinition}{\begin{definition}\begin{rm} }
\newcommand{\edefinition}{ \end{rm} \rule{1mm}{2mm}
\end{definition} }
\newcommand{\bproposition}{\begin{proposition} \begin{rm} }
\newcommand{\eproposition}{ \end{rm} \rule{1mm}{2mm}
\end{proposition} }
\newcommand{\bexample}{\begin{example} \begin{rm} }
\newcommand{\eexample}{ \end{rm} \rule{1mm}{2mm}
\end{example} }
\newcommand{\basm}{\begin{assumption} \begin{rm}}
\newcommand{\easm}{\end{rm} %\rule{1mm}{2mm}
\end{assumption}}

\newcommand{\real}{\mathds{R}}

\newcommand{\CC}{\mathcal{C}}
\newcommand{\DD}{\mathcal{D}}
\newcommand{\EE}{\mathcal{E}}

\newcommand{\GG}{\mathcal{G}}

\newcommand{\TT}{\mathcal{T}}

\newcommand{\LL}{\mathcal{L}}

\newcommand{\NN}{\mathcal{N}}
\newcommand{\VV}{\mathcal{V}}

\newcommand{\argmax}{\operatorname{argmax}}
\newcommand{\trace}{\operatorname{tr}}

\newtheorem{theorem}{\bf Theorem}[section]
\newtheorem{lemma}{\bf Lemma}[section]
\newtheorem{definition}{\bf Definition}[section]
\newtheorem{remark}{\bf Remark}[section]
\newtheorem{corollary}{\bf Corollary}[section]
\newtheorem{assumption}{\bf Assumption}[section]

\newcommand\oprocendsymbol{\hbox{$\blacksquare$}}
\newcommand\oprocend{\relax\ifmmode\else\unskip\hfill\fi\oprocendsymbol}

\date{}

\begin{document}
\setlength{\abovedisplayskip}{5pt} % remove align margin
\setlength{\belowdisplayskip}{5pt}
\setlength{\abovedisplayshortskip}{0pt}
\setlength{\belowdisplayshortskip}{0pt}

\begin{frontmatter}
\title{Attack-resilient Estimation of Switched Nonlinear Stochastic Cyber-Physical Systems}

%\thanks{Hunmin Kim and Minghui Zhu are with the School of Electrical Engineering and Computer Science, Pennsylvania State University, 201 Old Main, University Park, PA, 16802}

\author{Hunmin Kim$^*$}\ead{hunminkim3@gmail.com} \ and    % Add the
\author{Pinyao Guo$^{\dagger}$}\ead{pug132@ist.psu.edu} \ and
\author{Minghui Zhu$^*$}\ead{muz16@psu.edu} \ and              % e-mail address
\author{Peng Liu$^{\dagger}$}\ead{pliu@ist.psu.edu} 

\address{$^*$The School of Electrical Engineering and Computer Science,
Pennsylvania State University\\
%201 Old Main, University Park, PA, 16802\\
$^{\dagger}$The College of Information Sciences and Technology,
Pennsylvania State University\\
201 Old Main, University Park, PA, 16802\\
}

\thanks{This work is partially supported by the National Science Foundation (CNS-1505664).}
\maketitle

\begin{abstract}
This paper studies attack-resilient estimation of a class of switched nonlinear systems subject to stochastic noises. 
The systems are threatened by both of signal attacks and switching attacks.
%We consider two classes of attacks which are signal attacks and switching attacks.
The problem is formulated as the joint estimation of states, attack vectors and modes of hidden-mode switched systems.
We propose an estimation algorithm which is composed of a bank of state and attack vector estimators and a mode estimator.
The mode estimator selects the most likely mode based on modes' posterior probabilities induced by the discrepancies between obtained outputs and predicted outputs.
%We propose an estimation algorithm in which multiple estimators generate state and attack vector estimates for each modes, and a mode estimator chooses the most likely mode.
%which consists of multiple unknown input and state estimators with a mode estimator where each estimator generates estimates mode, and the mode estimator
We formally analyze the stability of estimation errors in probability for the proposed estimator associated with the true mode when the hidden mode is time-invariant but remains unknown.
%The stability and mean square estimation error bounds for the proposed estimator are formally analyzed for the true mode when the hidden mode is time-invariant but remains unknown.
For hidden-mode switched linear systems, we discuss a way to reduce computational complexity which originates from unknown signal attack locations.
Lastly, we present numerical simulations on
the IEEE 68-bus test system to show the estimator performance for time-varying modes with a regular mode set and a reduced mode set.
\end{abstract}

\end{frontmatter}
%\section{Comment} Equation~\eqref{CD001a} can be rewritten by \begin{align*} &\dot{x}(t) = f(x(t),u(t),K_s^{j(t)}d_a(t),w'(t),j(t),t),&&x(t) \in {\CC}^{j(t)}\nnum\\ &(x(t),j(t))^+=\Omega(x(t),j(t)),&&x(t) \in {\DD}^{j(t)}\nnum\\ &y_k = h(x_k,u_k,v_k',j_k,t_k)+K_H^{j_k} d'_{k,s}. \end{align*} Since $d_1$ should be reflected by $y_k$, it would be a part of $d'_{k,s}$. Since $d'_{k,s}$ is not reflected in the differential equation, $G_1=0$.
\section{Introduction}
Cyber-Physical Systems (CPS) are systems which integrate control systems with advanced technologies of sensing, computation and communication.
The integration leads to highly automated and collaborative applications such as autonomous vehicles, remote patient monitoring, and smart grid.
Due to such potentials, the significance on CPS is emphasized by the President's Council of Advisors on Science and Technology at 2010~\cite{holdren2010report}.
%The realizations of the CPS are accelerated by the developments of highly capable sensors, computing devices, and wireless communications~\cite{lee2010cyber}.
%However, there still remain many challenges such as resilience, privacy, interoperability, etc where some details can be found in~\cite{baheti2011cyber,cardenas2009challenges}.

%Secure control/attack detection~\cite{cardenas2009challenges}
%~\cite{cardenas2008research}
Security is of vital importance for CPS. Especially, because of the couplings between the cyber layer and physical layer, CPS bear vulnerabilities to cyberattacks which may cause irreparable damage to the physical layer~\cite{cardenas2008secure}.
%The increasingly internetworked CPS bear more vulnerabilities threatening reliable operations.
For example, a natural gas flow control system in Russia was temporarily seized in 2000
and a sewage control system in Australia was attacked in the same year~\cite{slay2007lessons}.
%For example, natural gas flow control system was temporary sized in Russia, and sewage control system was under attack in Australia~\cite{slay2007lessons} causing a flood nearby area.
According to early studies on CPS security, the types of possible attacks on CPS can be categorized into signal attacks and switching attacks.
Signal attacks include sensor attacks which tamper with sensor readings and actuator attacks which tamper with control commands. 
%Signal attacks target sensors/actuators or their communication lines. False data injection attacks~\cite{cardenas2008secure,mo2010false} or data replay attacks~\cite{mo2009secure} are examples.
While signal attacks modify the magnitudes or timings~\cite{zhu2013distributed,zhu2014performance} of signals, switching attacks alter system structures~\cite{weimer2012distributed,yong2015resilient}.
The attacks can be launched via communication jamming and malware; e.g., Trojan.
%Replay attacks tamper with the timing of message to be delivered~\cite{zhu2012resilient,zhu2013distributed,zhu2014performance}.

\textbf{Literature review.}
One topic in CPS security focuses on identifying fundamental limitations on the detectability of attacks.
Paper~\cite{liu2011false} investigates the possibilities that attacks can bypass detection algorithms for power systems.
Undetectable attack spaces on linear systems are characterized in~\cite{pasqualetti2013attack}.
A number of attack detectors against signal attacks are designed. In particular, attack detection problems for deterministic linear systems are formalized into $\ell_0$/$\ell_\infty$ optimization problems~\cite{fawzi2014secure,pasqualetti2013attack}.
Due to non-convexity, the problems are NP-hard~\cite{pasqualetti2013attack}.
To address the computational challenges, paper~\cite{fawzi2014secure} proposes convex relaxations of the optimization problems.
It is analyzed in~\cite{pajic2014robustness} that a state estimator based on $\ell_0$ optimization is robust with respect to modeling errors caused by sampling, computation/actuation jitter, and synchronization.
Paper~\cite{mo2010false} uses the Kalman filter to conduct attack-resilient state estimation in the presence of stochastic noises.
Paper~\cite{shoukry2016smt} proposes a multi-modal Luenberger observer whose memory usage increases linearly with the number of states and outputs.
Papers~\cite{liu2011false,mo2010false,pajic2014robustness,shoukry2016smt} consider sensor attacks, and
papers~\cite{fawzi2014secure,pasqualetti2013attack} consider both sensor and actuator attacks.
Paper~\cite{yong2015resilient} designs an attack-resilient estimator for stochastic linear systems in the presence of sensor attacks, actuator attacks, and switching attacks.  
All aforementioned papers focus on linear systems.

%There are few scalable resilient state estimation algorithms due to exponentially increasing computational complexity, except~\cite{shoukry2016smt}. We, in this paper, discuss a way to reduce computational complexity for switched linear time-varying systems, under both actuator attacks, and sensor attacks.

Our attack-resilient estimator design method is based on 
simultaneous unknown Input and State Estimation (ISE).
Early research of this area focuses on state estimation without estimating unknown inputs~\cite{darouach1997unbiased,hou1998optimal}.
Unbiased and minimum variance unknown input and state estimators are designed for linear systems without direct feedthrough matrix~\cite{hsieh2000robust} and with full-column rank direct feedthrough matrix~\cite{gillijns2007unbiased,yong2013simultaneous},
and with rank-deficient direct feedthrough matrix~\cite{SY-MZ-EF:AUTOMATICA16}.
Noticeably, this set of papers is restricted to linear systems.
%To our best knowledge, there is no systematic study on the ISE for nonlinear systems.
%Those unknown input and state estimation techniques are developed on linear systems. Recently,~\cite{yong2015resilient,yong2016simultaneous} design attack-resilient state, attack vector, and mode  estimator which runs the ISE in parallel and compare the estimation results. We propose a nonlinear extension of those in~\cite{yong2015resilient,yong2016simultaneous}. To achieve this, we extend the ISE to Nonlinear unknown Input and State Estimator (NISE) and we formally analyze the stability and estimation error bounds of the NISE to guarantee the attack-resilient estimation performance.

\textbf{Contributions.}
In this paper, we consider a class of switched nonlinear stochastic systems under signal attacks and switching attacks. We formulate the attack-resilient estimation problem as the simultaneous estimation of states, attack vectors and hidden modes.
The proposed algorithm associates an estimator to each mode and the estimators share the same structure.
Each estimator recursively produces the estimates of states and attack vectors.
The mode estimator selects the most likely mode based on modes' posteriori probabilities found by the discrepancies between obtained outputs and predicted outputs.
For the special case where the hidden mode is fixed but remains unknown,
it is shown that the estimation errors of states and attack vectors of the true mode satisfy Practically Exponentially Stable in probability (PESp) like properties.
In addition, we discuss a way to reduce computational complexity by reducing the number of modes for switched linear stochastic systems, which maintains the minimal number of modes to achieve the same detection capabilities as the power set. 
On the IEEE 68-bus test system, numerical simulations are conducted to show the effectiveness of the proposed algorithm on time-varying modes with a regular mode set and a reduced mode set.
Towards our best knowledge, this is the first time to systematically study unknown input, state, and mode estimation of switched nonlinear stochastic systems.

This paper is enriched from preliminary version~\cite{HK-MZ:ACC17} and includes a set of new results. In particular, this paper discusses how to reduce computational complexity caused by unknown signal attack locations.
% a large number of modes.
An additional numerical simulation is conduct to show the estimator performance for reduced mode sets. 
Lastly, this paper contains all the proofs which were omitted in~\cite{HK-MZ:ACC17}.

%Due to the first contribution, one can apply the ISE to a more general set of system. Because of the second contribution, we are able to apply well-studied properties of the EKF to the NISE problems such as stabilities, and convergences~\cite{kluge2010stochastic,xu2008stochastic}.
% Connecting two estimation problems is not trivial. In the current ISE, the state error update law seems to similar to that of the Kalman filter. However, there are two differences; in the ISE, the optimal gain should be found by constrained optimization problem and process error and output error are correlated. We address the issues by dividing the output to two(or three) parts by conducting a coordinate transformation basing on singular value decomposition. In doing so, we could subtract unknown input corrupted parts and uncorrupted parts. Each divided outputs have uncorrelated errors and they are used to estimate different state or inputs so that estimate errors are uncorrelated each other.

\textbf{Paper organization.} A motivating example of CPS model and attack model is introduced in Section~\ref{sec:motivating_ex}. 
Section~\ref{sec:power_exm} introduces system model, attack model and defender's knowledge. Moreover, the state, attack vector, and mode estimation problem is formulated in the same section.
We propose Nonlinear unknown Input, State and Mode Estimator (NISME) to solve the estimation problem in Section~\ref{sec:esti_deg}.
The stability of the proposed estimator is formally analyzed in Section~\ref{sec:analysis0}.
Section~\ref{sec:DiscR} discusses a way to reduce computational complexity caused by unknown signal attack locations, for switched stochastic linear systems.
Section~\ref{sec:nu_sim} offers numerical simulations on the IEEE 68-bus test system to show the performance of the proposed NISME with a regular mode set and a reduced mode set.
We present the derivations and proofs of the proposed estimator in Section~\ref{sec:proofs}.

\textbf{Notations.}
%We consider two types of unknown inputs. Those are direct feed through $d_{1,k}$ (or $d_1(t)$), and indirect feed through $d_{2,k}$ (or $d_2(t)$). Direct feed through input can be measured by output $y_k$ (or $y(t)$) but indirect feed through does not. 
Given a vector $a_{k}$, we use $\hat{a}_{k}$ and $\tilde{a}_{k}$ to denote an estimate of $a_{k}$ and induced estimation error $\tilde{a}_{k}=a_{k}-\hat{a}_{k}$, respectively. Its error covariance is defined by $P_k^{a} \triangleq {\mathbb E}[\tilde{a}_{k}\tilde{a}_{k}^T]$, and cross error covariance with $b_k$ is $P_k^{ab} \triangleq {\mathbb E}[\tilde{a}_{k}\tilde{b}_{k}^T]=(P_k^{ba})^T$.
%The number of $k-$combinations of a $n$-element set is denoted by $C_n^k$.

%Of $n$ elements, the number of $k$-combination is $C_n^k= \frac{n!}{k!(n-k)!}$. We define the set of all $k-$combinations of set $S$ as $\dbinom{S}{k}$.
We use the following definition for filter stability of nonlinear systems.
\begin{definition} 
Stochastic process $x(t)$ is said to be Practically Exponentially Stable in probability (PESp) if for any $\gamma \in (0,1)$, there exist positive constants $\alpha$, $b$, $c$, and $\delta$ such that, for any $\|x(0)\|\leq\delta$, the following holds for all $t\geq 0$:
\begin{align*}
P(\|x(t)\| < \alpha e^{-b t}\|x(0)\|+c) \geq 1-\gamma.
\end{align*}
\end{definition}
PESp is a special case of stochastic input-to-state stability~\cite{liu2008notion} when input is absent and class ${\mathcal K}{\mathcal L}$ function is exponential in $t$ and linear in $\|x(0)\|$. In addition, PESp is also extended from global asymptotic stability in probability (Definition 3.1 in~\cite{krstic1998stabilization}). Notice that the stability notions in~\cite{krstic1998stabilization,liu2008notion} are global and PESp is local.

%\begin{definition} \cite{agniel1971almost} Stochastic process $x(t)$ is said to be bounded with probability one, if it holds with probability one that $\sup_{t}\|x(t)\| < \infty$. \end{definition}
As for linear systems, one of the sufficient conditions for filter stability is uniform observability.
\begin{definition}
\cite{anderson1981detectability}
%\cite{reif1999stochastic}
The pair $(C_k,A_k)$ is uniformly observable if and only if there exist positive constants $a,b,l$, for all $k \geq 0$, such that, for all $k \geq 0$,
$a I \leq {\mathcal M}_{k+l,k} \leq b I$ where ${\mathcal M}_{k+l,k} \triangleq \sum_{i = k}^{k+l} \Phi_{i,k}C_iC_i^T \Phi_{i,k}^T$ is the observability gramian and $\Phi_{k_1,k_0}$ is the state transition matrix.
\end{definition}
Uniform observability reduces to observability if the linear system is time-invariant.
%\begin{definition} \cite{yong2016simultaneous} For any initial mode probability, assume that the geometric mean of the true mode probability asymptotically converges to 1. Then a filter is called mean consistency. \end{definition}
\section{Motivating example}\label{sec:motivating_ex}
%To motivate the mathematical description of the system model and possible attacks on the CPS, we bring the power transmission system and its possible vulnerabilities.
A power network is represented by undirected graph $({\VV},{\EE})$ with the set of buses ${\VV} \triangleq \{1,\cdots,N\}$ and the set of transmission lines ${\EE} \subseteq {\VV} \times {\VV}$.
The set of neighboring buses of $i \in {\VV}$ is ${\mathcal S}_i \triangleq \{l \in {\VV} \setminus \{i\} | (i,l) \in {\EE}\}$.
Each bus is either a generator bus $i \in {\GG}$, or a load bus $i \in {\LL}$. 
The dynamic of bus $i$ with attacks is described as the following switched nonlinear system:
\begin{align}
\dot{\theta}_i(t) &= f_i(t)+w_{1,i}(t)\nnum\\
\dot{f}_i(t) &= -\frac{1}{m_i}\big(D_i f_i(t) + \sum_{l\in {\mathcal S}_i} P_{il}^{j_{il}(t)}(t) \nnum\\
&-(P_{M_i}(t) +d_{a,i}(t))+P_{L_i}(t)\big)+w_{2,i}(t)
%\nnum\\ \dot{P}_{M_i}(t) &= -\frac{1}{T_{{CH}_i}}\big(P_{M_i}(t) - P_{v_i}(t)\big)\nnum\\ \dot{P}_{v_i}(t) &= -\frac{1}{T_{G_i}}\big(P_{v_i}(t) + \frac{1}{R_i} w_i(t) - P_{{ref}_i}(t)\big)
\label{Po_model}
\end{align}
with the output model
\begin{align}
y_{i,k} &= [P_{elec_i,k},\theta_{i,k},f_{i,k}]^T+d_{s,i,k}+v_{i,k}
\label{Po_model_output}
\end{align}
adopted from Chapter 9 in~\cite{Wood:1996} adding a phase angle measurement as~\cite{de2010synchronized,phadke1993synchronized}.
%States $\theta_i(t)$, $f_i(t)$, and variables $P_{L_i}(t)$, $P_{M_i}(t)$, $P_{il}^{j_{il}(t)}(t)$ are deviations from their nominal values.
System states $\theta_i(t)$, $f_i(t)$ are phase angle and angular frequency, respectively. 
Mode index $j_{il}(t) \in \{0,1\}$ represents on/off of the power line connection between buses $i$ and $l$; i.e., power flow is $P_{il}^{1}(t) =-P_{li}^{1}(t)= t_{il}\sin(\theta_i(t)-\theta_l(t))$, and $P_{il}^{0}(t) =-P_{li}^{0}(t)=0$.
The values $P_{L_i}(t)$, and $P_{elec_i}(t)= P_{L_i}(t)+D_if_i(t)$ denote power demand, and electrical power output, respectively. 
Since power demand $P_{L_i}(t)$ can be obtained by many load forecasting methods~\cite{alfares2002electric,hippert2001neural}, it is assumed to be known.

Mechanical power $P_{M_i}(t)$ is the control input for $i \in \GG$ and is assumed to be zero at load bus $i \in \LL$. %Since the power demand is known, $\bar{P}_{M_i}(t) = P_{M_i}(t)-P_{L_i}(t)$ can be used as the input of the system after an input transformation.
Power demand can be divided into elastic demand $P_{L_i}^E(t)$ and inelastic demand $P_{L_i}^{IE}(t)$ as shown in~\cite{Alvarado.Meng.ea:01}; i.e., $P_{L_i}(t)=P_{L_i}^E(t)+P_{L_i}^{IE}(t)$.
Elastic demand $P_{L_i}^E(t)$ can be controlled via power pricing. Since we assume that the current load is known, we simplify that load bus $i \in \LL$ uses $P_{L_i}(t)$ as load controller.

The measurements are sampled at discrete instants due to hardware constraints. We use subscript $k \in {\mathbb Z}_{\geq0}$ to denote an instantaneous value at the discrete sampling time $t_k$; e.g., $f_i(t_k) =f_{i,k}$.
%In the output, we can modify output $P_{elec_i,k}$ into $\frac{P_{elec_i,k}-P_{L_i,k}}{D_i}$ because power demand $P_{L_i,k}$ and load damping constant $D_i$ are known.

An attacker is assumed to be able to modify the sensor measurements, control commands, and trigger the power flow line switches. The possible attacks are modeled as vectors $d_{s,i,k} \in {\mathbb R}$, $d_{a,i}(t) \in {\mathbb R}$, and hidden mode switch $j_{il}(t)$ which represent sensor attacks~\cite{mo2010false,pasqualetti2013attack}, actuator attacks~\cite{fawzi2014secure,pasqualetti2013attack,zhu2014performance}, and circuit breaking/switching attacks~\cite{weimer2012distributed,yong2015resilient}, respectively.

%Sensor attack $d_{s,i,k}$ is a type of false data injection attack~\cite{mo2010false,pasqualetti2013attack} or data replay attack~\cite{mo2009secure} where the attack target is the sensor data. Actuator attack $d_{a,i}(t)$ represents another type of false data injection attack~\cite{fawzi2014secure,pasqualetti2013attack} or data replay attack~\cite{mo2009secure,zhu2014performance} where the targeted signal is the control input. Power line switch attack $j_{il,k}$ triggers circuit breakers, breaks the power line physically, or turns on/off power line switch to modify the power flow~\cite{weimer2012distributed,yong2015resilient}; i.e., $P_{il}^{j_{il,k}=0}=0$ and $P_{il}^{j_{il,k}=1}=t_{il}\sin(\theta_i(t)-\theta_l(t))$. The attack is modeled as a hidden-mode change because the current mode and its transition are unknown to the defender.

\section{Problem formulation}\label{sec:power_exm}
\textbf{System model.} Consider the hidden-mode nonlinear stochastic system
\begin{align}
&\dot{x}(t) = f'(x(t),u(t)+d_a(t),w'(t),j(t),t), &&x(t) \in {\CC}^{j(t)}\nnum\\
&(x(t),j(t))^+=\Omega'(x(t),j(t)), &&x(t) \in {\DD}^{j(t)}\nnum\\
&y_k = h(x_k,u_k+d_{a,k},v_k',j_k,t_k)+d_{s,k}
\label{CD001z}
\end{align}
where $x(t) \in {\real}^{n}$, $y_k \in {\real}^{m}$, $u(t) \in {\real}^{s}$, and $j(t) \in {\mathbb M}^I$ are state, output, input, and hidden-mode, respectively.
We use subscript $k \in {\mathbb Z}_{\geq0}$ to denote an instantaneous value at the discrete sampling time $t_k$.
Vectors $d_a(t) \in {\real}^{s}$, and $d_{s,k} \in {\real}^{m}$ are actuator attack vector and sensor attack vector, respectively.
Sets ${\CC}^{j(t)}$, ${\DD}^{j(t)} \subseteq {\real}^{n}$ denote flow set and jump set, respectively, and $\Omega$ is a mode transition function. For each mode, process noise $w'(t)\in {\real}^{s_1}$ and measurement noise $v_k' \in {\real}^{s_2}$ are uncorrelated with each other. %Their covariance matrices are given by $Q_k'^{j_k} \triangleq {\mathbb E}[v_k'v_k'^T]$, and $R_k'^{j_k} \triangleq {\mathbb E}[w_k'w_k'^T]$.
The system is a continuous-discrete system because, while the physical dynamic evolves in continuous time, sensor measurements are obtained at their corresponding sampling instants due to hardware constraints.
%We define such instant as $t_k$ for $k \in {\mathbb Z}_{\geq0}$, and a uniform sampling period as $\epsilon = t_{k}-t_{k-1}$.
We define a uniform sampling period as $\epsilon = t_{k}-t_{k-1}$.
It is assumed that the system~\eqref{CD001z} has a unique solution. One of the sufficient condition for the unique solution is weak one-sided local Lipschitz condition on function $f'(\cdot)$ in the open time interval of each mode duration~\cite{von2010existence} and other conditions can be found in the references therein.
The system model~\eqref{CD001z} includes the power system model~\eqref{Po_model} with~\eqref{Po_model_output} as a special case.

\textbf{Attack model.}
Signal attacks are comprised of signal magnitude attacks (i.e., the attacker injects attack signals), and signal location attacks (i.e., the attacker chooses targeted sensors/actuators).
Signal attacks are modeled by $d_a(t)$ and $d_{s,k}$ where zero values indicate that the corresponding actuators and sensors are free of attacks and non-zero values represent attack magnitudes.
Switching attacks change system modes following $\Omega'$.
%We model signal attacks as actuator attack vector $d_a(t)$ and sensor attack vector $d_{s,k}$, where their magnitude represents magnitude attacks, and a combination of non-zero elements in the vectors represents signal location attacks. Switching attacks change the mode of the system operations, and measurements.
%We denote mode index $j'(t)$ to model the system operational changes. and measurement model changes. Moreover the modes might change over time, and we make no assumption on attacks.

\textbf{Knowledge of the defender.} 
The defender is unaware of which actuators/sensors are under attacks and what the current mode is.
The defender knows dynamic system model and output model~\eqref{CD001a} for each mode but not the mode transition function $\Omega'$. 
Mode set ${\mathbb M}^I$ is also known to the defender.
%If the defender does not have any knowledge on attack location ${\mathbb M}^A$, the combination of all $n+m$ possible locations are considered; i.e., $2^{n+m}$ in total.
%it can be constructed by setting $2^{n+m}$ numbers of modes because the number of locations is only $n+m$.
%The attack vector might be a predetermined signal or a random signal with an unknown distribution.
The attack vectors $d_a(t)$, $d_{s,k}$, mode $j(t)$ and its transitions are inaccessible to the defender.
Noise vectors $w'(t)$, $v_k'$ are unknown but their auto covariance matrices are known.

%Given the system~\eqref{CD001a} with attacks, the problem for the defender is to estimate the state resiliently as well as the attack vector, given a tuples of system functions/matrices, control input, sensor outputs, and previous step estimates. In other words, the defender designs estimator ${\TT} (\Lambda_k) = [\hat{x}_k^T,\hat{d}_k^T,P_{k}^x]^T$ to generate state and attack estimates and covariance matrix where the estimator input is given by $\Lambda_k = (f(\cdot,\cdot,t),h(\cdot,\cdot,t_k),H_k,u(t),y_{k},\hat{x}_{k-1}, \hat{d}_{k-1},P_{k-1}^x)$.

\textbf{Objective.} 
The defender aims to answer the following three questions:\\
(a) if any sensor or actuator is attacked;\\
(b) if so, which ones are attacked, and how much sensor readings and control commands are tampered with;\\
(c) what current system states are.\\
The above problem can be formulated as a joint estimation of states, attack vectors and modes of hidden-mode switched systems~\eqref{CD001z}.

%The problem for the defender is to design an estimator to generate state, attack vector, and mode estimates.

\section{Estimator design}\label{sec:esti_deg}
In order to reflect real world, system~\eqref{CD001z} models the attacks from the attacker's point of view and captures attack sources. In order to solve the estimation problem, we need to model the attacks from the defender's point of view and captures attack consequences. In particular, we rewrite system~\eqref{CD001z} as follows:
\begin{align}
&\dot{x}(t) = f(x(t),u(t),d(t),w'(t),j(t),t),&&x(t) \in {\CC}^{j(t)}\nnum\\
&(x(t),j(t))^+=\Omega(x(t),j(t)),&&x(t) \in {\DD}^{j(t)}\nnum\\
&y_k = h(x_k,u_k,v_k',j_k,t_k)+H^{j_k} d_k
\label{CD001a}
\end{align}
where $d(t) = [d_a^T(t),d_{s}'^T(t)]^T \in {\real}^{s+m}$, $d_{s,k}' = h(x_k,u_k+d_{a,k},v_k',j_k,t_k)-h(x_k,u_k,v_k',j_k,t_k)+d_{s,k}$, $f(x(t),u(t)$, $d(t),w'(t),j(t),t)$ $= f'(x(t),u(t)+S^{j(t)}d(t),w'(t)$, $j(t),t)$, $S^{j(t)} = [{K}_S^{j(t)}$, $0^{s \times m}] \in \{0,1\}^{s \times (s+m)}$ and $H^{j_k} = [0^{m \times s}, K_H^{j_k}] \in \{0,1\}^{m \times (s+m)}$.
The defender models the signal location attacks as mode $j(t)$ of diagonal matrix
\begin{align*}
K^j \triangleq 
\left[
\begin{array}{c}
S^{j}\\
H^{j}\\
\end{array}
\right] =
\left[
\begin{array}{cc}
K_S^{j}& 0^{s \times m}\\
0^{m \times s}& K_H^{j}\\
\end{array}
\right] \in \{0,1\}^{(s+m) \times (s+m)}
\end{align*}
where ${K}^{j}(i,i)=1$ if mode $j$ assumes that the $i^{th}$ location is under attack; otherwise, ${K}^{j}(i,i)=0$.
Thus, $j(t) \in {\mathbb M} = {\mathbb M^A} \times {\mathbb M^I}$ stands for the both signal location attacks ${\mathbb M^A}$ and switching attacks ${\mathbb M^I}$.

\begin{remark}
For the sake of generality, we will consider arbitrary $K^j \in R^{(s+m) \times (s+m)}$ in
the remaining of this section and Sections~\ref{sec:analysis0}.
\oprocend
\end{remark}

\begin{figure}[h]
  \centering
  \includegraphics[width = \linewidth]{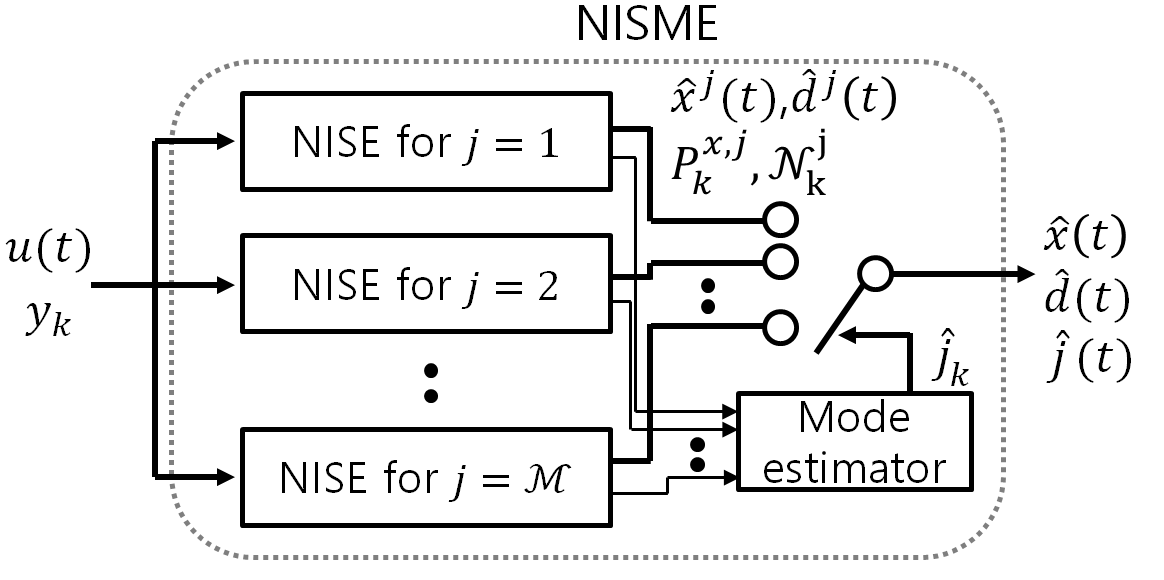}
  \caption{Scheme of NISME}\label{NISME_F}
\end{figure}

To solve the problem, we propose Nonlinear unknown Input, State and Mode Estimator (NISME).
The NISME consists of a bank of Nonlinear unknown Input and State estimators (NISE) and a mode estimator
%to elect the most likely one among the estimates generated by the NISEs
as shown in Figure~\ref{NISME_F}.
Each NISE is associated with a particular mode and recursively estimates the states and attack vectors under the fixed mode. The mode estimator calculates the posteriori probabilities of the modes by observing output discrepancies from predicted outputs, and chooses the most likely one. 
Lastly, the NISME outputs the estimates of the states and the attack vectors of the selected mode.

We first introduce some preliminaries for the NISE in Section~\ref{sys_model}.
The NISME is presented in Section~\ref{sec:Algo0}.

\subsection{Preliminaries}\label{sys_model}
In this section, we introduce an output decomposition used in the NISE. 
Since each NISE is associated with a particular mode, we omit the mode index $j(t)$ for notational simplicity.
%We keep any fixed mode $j_k = j \in {\mathbb M}$ for $\forall k$ for the filter derivation and analysis.
%For notational simplicity, we omit the mode notation $j_k$ if it is clear.

We first discretize and linearize system~\eqref{CD001a} as follows with constant sampling period $\epsilon$:
%We need a discretized system model to estimate attack vector $d_k$. The discretized and linearized models are, with constant sampling period $\epsilon=t_{k+1}-t_k$,
\begin{align}
%x(t) & \simeq A(t)x(t)+B(t)u(t)+G(t)d(t)+w(t)\nnum\\
x_{k+1} &= x_k + \epsilon f(x_k,u_k,d_k,w_k',t_k)+\epsilon \rho_{k} \nnum\\
&\simeq x_k+\epsilon (A_k x_k + B_ku_k+ G_kd_k+\rho_{k}+w_k)\nnum\\
y_k&\simeq C_kx_k+D_ku_k+H d_k+v_k
\label{E001}
\end{align}
where $\epsilon \rho_k \triangleq \int_{t_{k}}^{t_{k+1}}f(x(\tau),u(\tau),d(\tau),w'(\tau),\tau) d \tau -\epsilon f(x_k,u_k,d_k,w_k',t_k)$ refers to discretization error, and $w_k=J_{k} w_k'$, $v_k=E_k v_k'$,
\begin{align*}
&A_k \triangleq \frac{\partial f_k}{\partial x}\big|_{\hat{x}_{k|k},u_{k},\hat{d}_k,0,t_k},
&&B_k \triangleq \frac{\partial f_k}{\partial u}\big|_{\hat{x}_{k|k},u_{k},\hat{d}_k,0,t_k},\nnum\\
&G_k \triangleq \frac{\partial f_k}{\partial d}\big|_{\hat{x}_{k|k},u_{k},\hat{d}_{k},0,t_k},
&&J_{k} \triangleq \frac{\partial f_k}{\partial w'}\big|_{\hat{x}_{k|k},u_{k},\hat{d}_k,0,t_k},\nnum\\
&C_k \triangleq \frac{\partial h_k}{\partial x}\big|_{\hat{x}_{k|k-1},u_{k},0,t_k},
&&D_k \triangleq \frac{\partial h_k}{\partial u}\big|_{\hat{x}_{k|k-1},u_{k},0,t_k},\nnum\\
&E_{k}  \triangleq \frac{\partial h_k}{\partial v'}\big|_{\hat{x}_{k|k-1},u_{k},0,t_k}.
\end{align*}
We define the autocovariance matrices for noise vectors as ${\mathbb{E}}[w_kw_k^T]=Q_k \geq 0$, and ${\mathbb{E}}[v_kv_k^T]=R_k>0$.

Now we introduce two coordinate transformations. The first one is based on the singular value decomposition
\begin{align*}
&H=
\left[
\begin{array}{cc}
U_{1}&U_{2}\\
\end{array}
\right]
\left[
\begin{array}{cc}
{\Sigma}&0\\
0&0\\
\end{array}
\right]
\left[
\begin{array}{c}
{V}_{1}^T\\
{V}_{2}^T\\
\end{array}
\right]
\end{align*}
where $\Sigma$ is a full rank diagonal matrix. The first coordinate transformation $T_{k}$ is defined by
\begin{align}
T_{k}&=
\left[
\begin{array}{c}
T_{1,k}\\T_{2}
\end{array}
\right]=
\left[
\begin{array}{cc}
I& - U_{1}^T R_k U_{2} (U_{2}^TR_kU_{2})^{-1}\\
0&I
\end{array}
\right]
\left[
\begin{array}{c}
U_{1}^T\\
U_{2}^T\\
\end{array}
\right].
\label{P202.1}
\end{align}
Likewise, the singular value decomposition
\begin{align*}
T_{2}C_{k}G_{k-1}V_{2}=
\left[
\begin{array}{cc}
\bar{U}_{1,k}&\bar{U}_{2,k}\\
\end{array}
\right]
\left[
\begin{array}{cc}
\bar{\Sigma}_{k}&0\\
0&0\\
\end{array}
\right]
\left[
\begin{array}{c}
\bar{V}_{1,k}^T\\
\bar{V}_{2,k}^T\\
\end{array}
\right]
\end{align*}
with full-rank diagonal matrix $\bar{\Sigma}_k$ induces the second coordinate transformation $\bar{T}_{k}=[\bar{T}_{1,k}^T, \bar{T}_{2,k}^T]^T$
\begin{align}
\bar{T}_{k}&=
\left[
\begin{array}{cc}
I& - \bar{U}_{1,k}^T \bar{R}_k \bar{U}_{2,k} (\bar{U}_{2,k}^T\bar{R}_k\bar{U}_{2,k})^{-1}\\
0&I
\end{array}
\right]
\left[
\begin{array}{c}
\bar{U}_{1,k}^T\\
\bar{U}_{2,k}^T\\
\end{array}
\right]
\label{P202.2}
\end{align}
where $\bar{R}_k \triangleq T_{2}R_kT_{2}^T$.
From coordinate transformations~\eqref{P202.1} and~\eqref{P202.2}, the output $y_k$ in~\eqref{E001} can be decomposed as follows:
\begin{align}
z_{1,k}&=T_{1,k}y_k \simeq C_{1,k} x_k + D_{1,k}u_k + H_{1} d_{1,k} + v_{1,k}\nnum\\
z_{2,k}&=\bar{T}_{1,k}T_{2}y_k \simeq C_{2,k}x_k+D_{2,k}u_k+v_{2,k}\nnum\\
&\simeq C_{2,k}(x_{k-1}+\epsilon (A_{k-1} x_{k-1} + B_{k-1}u_{k-1}\nnum\\
&+ G_{k-1}d_{k-1}+w_{k-1}))+ D_{2,k}u_k+ v_{2,k}\nnum\\
z_{3,k}&=\bar{T}_{2,k}T_{2}y_k \simeq C_{3,k}x_k+ D_{3,k}u_k+ v_{3,k}
\label{CD002}
\end{align}
where $H_{1}=\Sigma$ and $C_{2,k}G_{k-1}V_{2}=\bar{\Sigma}_k$.
Attack vector $d_k$ is decomposed into a sum of $d_{1,k} \triangleq V_{1}^Td_k$ and $d_{2,k} \triangleq V_{2}^Td_k$ where they are orthogonal to each other.
Note that $d_{1,k}$ and $d_{2,k}$ are different from $d_{s,k}'$ and $d_{a,k}$, and introduced for the purpose of analysis.
In this case, it holds that $G_k d_k= G_{1,k}d_{1,k}+G_{2,k}d_{2,k}$ with $G_{1,k} \triangleq G_kV_{1}$ and $G_{2,k} \triangleq G_kV_{2}$.
Output $z_{1,k}$ is the portion of $y_k$ which is attacked at $k$; i.e., $z_{1,k}$ includes $d_{1,k}$ in~\eqref{CD002}.
Outputs $z_{2,k}$ and $z_{3,k}$ are the portions of $y_k$ and are free of attacks at $k$, where output $z_{2,k}$ reflects $d_{2,k}$ indirectly because $C_{2,k}G_{k-1}d_{k-1}=\bar{\Sigma}$.
Thus, decomposed outputs $z_{1,k}$, $z_{2,k}$, and $z_{3,k}$ are used to estimate $d_{1,k}$, $d_{2,k-1}$, and $x_k$, respectively.

Because $d_{2,k-1}$ is not measured by $y_{k-1}$, output $z_{2,k}=C_{2,k}G_{2,k-1}d_{2,k-1}+...$ in~\eqref{CD002} is instead used to estimate $d_{2,k-1}$; i.e., matrices $C_{2,k}$ and $G_{2,k-1}$ must be known to estimate attack vector $d_{2,k-1}$.
However, in~\eqref{E001}, matrix $G_{2,k-1}$ is obtained by linearizing $f(\cdot)$ using $\hat{d}_{2,k-1}$, and matrix $C_k$ is obtained by linearizing $h(\cdot)$ using $\hat{x}_{k|k-1}$, where these linearizations cannot be done without knowing $\hat{d}_{2,k-1}$.
Thus, we have the following assumption.
\begin{assumption}
Dynamic system model~\eqref{CD001a} can be expressed as
\begin{align}
\dot{x}(t) &= f(x(t),u(t),d_1(t),w'(t),t)+G_2(t)d_2(t)\nnum\\
z_{1,k}&=T_{1,k}y_k= h_1(x_k,u_k,v_{1,k}',t_k) + H_{1} d_{1,k}\nnum\\
z_{2,k}&=\bar{T}_{1,k}T_{2,k}y_k= C_{2,k}x_k + h_2(u_k,v_{2,k}',t_k)\nnum\\
z_{3,k}&=\bar{T}_{2,k}T_{2,k}y_k = h_3(x_k,u_k,v_{3,k}',t_k).
\label{E003}
\end{align}
\label{CD_asm_m}
\end{assumption}
With Assumption~\ref{CD_asm_m}, the dynamic system~\eqref{E001} becomes
\begin{align}
%x(t) & \simeq A(t)x(t)+B(t)u(t)+G(t)d(t)+w(t)\nnum\\
&x_{k+1} = x_k+\epsilon f(x_k,u_k,d_{1,k},w_k',t_k)+\epsilon G_{2,k}d_{2,k}+\epsilon \rho_{k}\nnum\\
&\simeq x_k +\epsilon (A_k x_k + B_ku_k+ G_{1,k}d_{1,k}+G_{2,k}d_{2,k}+\rho_{k}\nnum\\
&+w_k)
\label{DD001}
\end{align}
where matrices $A_k$, $B_k$, $G_{1,k}$, and $J_{k}$ can be obtained before having an estimate for $d_{2,k}$. 
Output equation~\eqref{E003} is linearized into~\eqref{CD002} where noises
$v_{1,k}=E_{1,k} v_{1,k}'$, $v_{2,k}=E_{2,k} v_{2,k}'$, $v_{3,k}=E_{3,k} v_{3,k}'$ are uncorrelated with each other and
\begin{align*}
&C_{1,k} \triangleq \frac{\partial h_{1,k}}{\partial x}\big|_{\hat{x}_{k|k-1},u_{k},0,t_k},
&&C_{3,k} \triangleq \frac{\partial h_{3,k}}{\partial x}\big|_{\hat{x}_{k|k-1},u_{k},0,t_k},\nnum\\
&D_{1,k} \triangleq \frac{\partial h_{1,k}}{\partial u}\big|_{\hat{x}_{k|k-1},u_{k},0,t_k},
&&D_{2,k} \triangleq \frac{\partial h_{2,k}}{\partial u}\big|_{u_{k},0,t_k},
\end{align*}
\begin{align*}
&D_{3,k} \triangleq \frac{\partial h_{3,k}}{\partial u}\big|_{\hat{x}_{k|k-1},u_{k},0,t_k},
&&E_{1,k}  \triangleq \frac{\partial h_{1,k}}{\partial v_1'}\big|_{\hat{x}_{k|k-1},u_{k},0,t_k},\nnum\\
&E_{2,k}  \triangleq \frac{\partial h_{2,k}}{\partial v_2'}\big|_{u_{k},0,t_k},
&&E_{3,k}  \triangleq \frac{\partial h_{3,k}}{\partial v_3'}\big|_{\hat{x}_{k|k-1},u_{k},0,t_k}.
\end{align*}

\subsection{Algorithm statement}\label{sec:Algo0}
\begin{figure}[h]
  \centering
  \includegraphics[width = \linewidth]{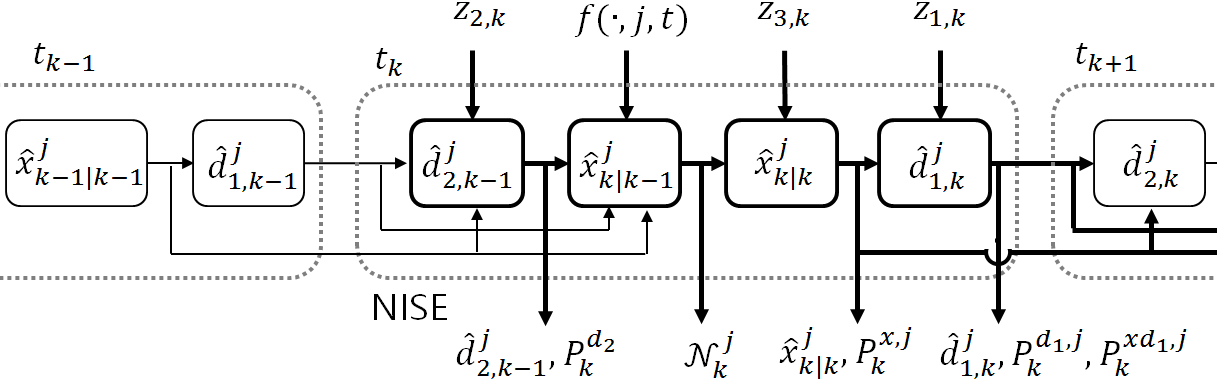}
  \caption{The recursive estimation scheme of the NISE.}\label{NISE_F}
\end{figure}
\begin{algorithm} \caption{NISE} \label{algo}
\begin{algorithmic}[1]
\REQUIRE $j,\hat{x}_{k-1|k-1}^j,\hat{d}_{1,k-1}^j$, $P_{k-1}^{x,j}$, $P_{k-1}^{d_1,j}$, $P_{k-1}^{xd_1,j}$, $y_k$, $u(t)$ for $t \in [t_{k-1},t_k]$;\\
$\triangleright$ \textbf{\textit{Attack vector $\textbf{d}_{2,k-1}^{j}$ estimation}}
\STATE $M_{2,k}^j=(\epsilon C_{2,k}^jG_{2,k-1}^j)^{-1}$ (or $M_{2,k}^j=0$ if $Rank(\bar{\Sigma}_k)=0$);
\STATE $\hat{d}_{2,k-1}^j = M_{2,k}^j (z_{2,k}^j-C_{2,k}^j (\hat{x}_{k-1|k-1}^j + \epsilon f(\hat{x}_{k-1|k-1}^j,$ $u_{k-1},\hat{d}_{1,k-1}^j,0,j,t_{k-1}))- h_2(u_k,0,j,t_k))$;
\STATE 
$P_{k-1}^{d_2,j}= M_{2,k}^jC_{2,k}^j(I+\epsilon A_{k-1}^j)P_{k-1}^j$ $(M_{2,k}^jC_{2,k}^j(I+\epsilon A_{k-1}^j))^T+\epsilon^2M_{2,k}^jC_{2,k}^jQ_{k-1}^j(M_{2,k}^j$ $C_{2,k}^j)^T+M_{2,k}^jR_{2,k}^j$ $(M_{2,k}^j)^T+\epsilon^2 C_{2,k}^j G_{1,k-1}^j P_{k-1}^{d_1,j}(C_{2,k}^j G_{1,k-1}^j)^T+\epsilon M_{2,k}^j C_{2,k}^j (I+\epsilon A_{k-1}^j)P_{k-1}^{x d_1,j}(C_{2,k}^jG_{1,k-1}^j)^T
+\epsilon C_{2,k}^jG_{1,k-1}^j P_{k-1}^{d_1 x,j} (M_{2,k}^j C_{2,k}^j (I+\epsilon A_{k-1}^j))^T$;\\
%$P_{k-1}^{d_2,j}=(\epsilon^2 (C_{2,k}^jG_{2,k-1}^j)^T(\tilde{R}_{2,k}^j)^{-1}C_{2,k}^jG_{2,k-1}^j)^{-1}$ where $\tilde{R}_{2,k}^j = C_{2,k}^j\tilde{P}_{k-1}^j(C_{2,k}^j)^T + R_{2,k}^j$, and $\tilde{P}_{k-1}^j=(I+\epsilon A_{k-1}^j) P_{k-1}^{x,j}(I+\epsilon A_{k-1}^j)^T- \epsilon(I+\epsilon A_{k-1}^j)P_{k-1}^{xd_1,j}(G_{1,k-1}^j)^T+ \epsilon G_{1,k-1}^jP_{k-1}^{d_1x,j}(I+\epsilon A_{k-1}^j)^T+ \epsilon^2 G_{1,k-1}^jP_{k-1}^{d_1,j}(G_{1,k-1}^j)^T+Q_{k-1}^j$;
$\triangleright$ \textbf{\textit{State prediction}}
\STATE $\dot{\hat{x}}^j(t) = f(\hat{x}^j(t),u(t),\hat{d}_{1,k-1}^j,0,j,t)+G_{2,k-1}^j\hat{d}_{2,k-1}^j$ with initial condition $\hat{x}_{k-1|k-1}^j$ for $t \in (t_{k-1},t_k]$ to have $\hat{x}_{k|k-1}^j$ at $t=t_k$;
\STATE $P_{k|k-1}^{x,j}=\bar{A}_{k-1}^jP_{k-1}^{x,j}(\bar{A}_{k-1}^j)^T+\bar{Q}_{k-1}^j$;\\
$\triangleright$ \textbf{\textit{State estimation}}
\STATE $L_k^j=P_{k|k-1}^{x,j}(C_{3,k}^j)^T (C_{3,k}^jP_{k|k-1}^{x,j}(C_{3,k}^j)^T+R_{3,k}^j)^{-1}$;
\STATE $\hat{x}_{k|k}^j = \hat{x}_{k|k-1}^j+L_{k}^j(z_{3,k}^j -h_3(\hat{x}_{k|k-1}^j,$ $u_k,0,j,t_k))$;
\STATE $P_{k}^{x,j} =(I-L_{k}^jC_{3,k}^j)P_{k|k-1}^{x,j}(I-L_{k}^jC_{3,k}^j)^T +L_k^jR_{3,k}^j(L_k^j)^T$;\\
$\triangleright$ \textbf{\textit{Attack vector $d_{1,k}^j$ estimation}}
\STATE $M_{1}^j =(H_{1}^j)^{-1}$ (or $M_{1}^j=0$ if $Rank(\Sigma)=0$);
\STATE $\hat{d}_{1,k}^j = M_{1}^j(z_{1,k}^j-h_1(\hat{x}_{k|k}^j,u_k,0,j,t_k))$;
\STATE 
$P_{k}^{d_1,j}=M_{1}^jC_{1,k}^jP_{k}^{x,j}(M_{1}^jC_{1,k}^j)^T+M_{1}^jR_{1,k}^j(M_{1}^j)^T$;
%$P_{k}^{d_1,j}=((H_{1,k}^j)^T(\tilde{R}_{1,k}^j)^{-1}H_{1,k}^j)^{-1}$
\STATE $P_{k}^{xd_1,j}=-P_{k}^{x,j}(M_{1}^jC_{1,k}^j)^T$;\\
%where $\tilde{R}_{1,k}^j=C_{1,k}^jP_{k}^{x,j}(C_{1,k}^j)^T+R_{1,k}^j$;
$\triangleright$ \textbf{\textit{The priori probability of the mode}}
\STATE $\nu_k^j = z_{3,k}^j -h_3(\hat{x}_{k|k-1}^j,u_k,0,j,t_k)$;
\STATE $\bar{P}_{k|k-1}^j = C_{3,k}^jP_{k|k-1}^j(C_{3,k}^j)^T+R_{3,k}^j$;
\STATE ${\NN}_k^j = \frac{1}{(2 \pi)^{n^j/2} |\bar{P}_{k|k-1}^j|^{1/2}}\exp(-\frac{(\nu_k^j)^T(\bar{P}_{k|k-1}^j)^{-1}\nu_k^j}{2})$;
\ENSURE $\hat{x}_{k|k}^j$, $\hat{d}_{1,k}^j$, $\hat{d}_{2,k-1}^j$, $P_{k}^{x,j}$, $P_{k-1}^{d_2,j}$, $P_{k}^{d_1,j}$, $P_{k}^{xd_1,j}$, ${\NN}_k^j$.
\end{algorithmic}
\end{algorithm}

Consider the NISE (Algorithm~\ref{algo}) as well as Figure~\ref{NISE_F} whose derivation is presented in Section~\ref{CD_filter} in details.
%We have the dynamic system~\eqref{CD001a} and transformed outputs~\eqref{CD002}. 
All the estimates of states and attack vectors are best linear unbiased estimates (BLUE); i.e., the estimator gains are chosen such that the estimates are unbiased and the norms of the error covariance matrices are minimized.
Since attack vector $d_{2,k-1}$ does not influence output $y_{k-1}$ directly, output $z_{2,k}^j$ is used to estimate attack vector $d_{2,k-1}$ (line 2) by using previous estimate of attack vector $d_{1,k-1}$. Error covariance matrix $P_{k-1}^{d_2,j}$ of attack vector estimate $\hat{d}_{2,k-1}^j$ is derived in line 3.
Applying the previous state and attack vector estimates to dynamic system~\eqref{CD001a}, the current state is predicted (line 4).
Error covariance matrix $P^{x,j}_{k|k-1}$ of the predicted state $\hat{x}_{k|k-1}^j$ is found in line 5 and the matrices in line 5 are defined by $\bar{Q}_{k-1}^j \triangleq {\mathbb E}[\bar{w}_{k-1}^j(\bar{w}_{k-1}^j)^T]$,
\begin{align}
\bar{w}&_{k-1}^j \triangleq \epsilon (I-\epsilon G_{2,k-1}^jM_{2,k}^jC_{2,k}^j)\nnum\\
&\times(w_{k-1}^j-G_{1,k-1}^jM_{1}^jv_{1,k-1}^j)-\epsilon G_{2,k-1}^jM_{2,k}^j v_{2,k}^j\nnum\\
\bar{A}&_{k-1}^j \triangleq (I-\epsilon G_{2,k-1}^jM_{2,k}^j C_{2,k}^j)\nnum\\
&\times(I+\epsilon A_{k-1}^j-\epsilon G_{1,k-1}^jM_{1}^jC_{1,k-1}^j).
\label{CD004.7}
\end{align}
We correct the predicted state using the measurement bias (line 7) between the measured output and the predicted output. Error covariance matrix $P_{k|k-1}^{x,j}$ of state prediction $\hat{x}_{k|k-1}^j$ is updated in line 8.
Attack vector $d_{1,k}$ is estimated from output $z_{1,k}^j$ (line 10). Error covariance matrix $P_{k-1}^{d_1,j}$ of attack vector estimate $\hat{d}_{1,k}^j$, and cross error covariance matrix $P_{k}^{xd_1,j}$ with $\hat{x}_{k|k}^j$ are found in lines 11-12.
Lastly, the NISE generates the priori probability ${\NN}_k^j$ (line 15) of the mode to find the most likely mode, where $n^{j} \triangleq Rank(\bar{P}_{k|k-1}^{j})$. For this purpose, the discrepancy between the measured output $z_{3,k}^j$ and the predicted output is used to validate the mode (line 13) because 
they should match if $j$ is the true mode.
Since the system is nonlinear, the discrepancy $\nu_k^j$ may not be Gaussian.
We approximate $\nu_k^j$ as a Gaussian random vector because it is a typical practice to approximate an unknown noise as a Gaussian distribution as~\cite{kotecha2003gaussian}. Moreover, $\nu_k^j$ is Gaussian when the system is linear 
and noises $w_k^j$ and $v_k^j$ are Gaussian.
Covariance matrix $\bar{P}_{k|k-1}^{j}$ of the discrepancy $\nu_k^j$ is found in line 14.

\begin{algorithm} \caption{NISME} \label{algo_mode}
\begin{algorithmic}[1]
\REQUIRE $\hat{x}_{0|0}^j={\mathbb E}[x_0]$, $P_0^{x,j}=P_0$, $\mu_0^j=\frac{1}{|{\mathbb M}|}$, $\hat{d}_{1,0}^j=(\Sigma^j)^{-1}(z_{1,0}^j-h_1(\hat{x}_{0|0}^j,u_0,0,j,0))$
for $\forall j \in {\mathbb M}$; Choose $0<\delta \ll \frac{1}{|{\mathbb M}|}$, $0<\alpha_1<1$, $0<\alpha_2<1$ (significance levels);
\FOR{$k=1:N$}
\STATE Read sensor output $y_k$, and control input $u(t)$ for $t \in [t_{k-1},t_{k}]$;
\FOR{$j \in {\mathbb M}$}
\STATE Run the NISE with input $(j,\hat{x}_{k-1|k-1}^j,\hat{d}_{1,k-1}^j,$ $P_{k-1}^{x,j}, P_{k-1}^{d_1,j}, P_{k-1}^{xd_1,j},y_k,u(t) {\rm \ for \ } t \in [t_{k-1},t_k])$ to generate output ($\hat{x}_{k|k}^j$, $\hat{d}_{1,k}^j$, $\hat{d}_{2,k-1}^j$, $P_{k}^{x,j}$, $P_{k-1}^{d_2,j}$, $P_{k}^{d_1,j}$, $P_{k}^{xd_1,j}$, ${\NN}_k^j$);
\ENDFOR \\
$\triangleright$ \textbf{\textit{Mode estimator}}
\FOR{$j \in {\mathbb M}$}
\STATE $\bar{\mu}_k^j = \max\{\frac{{\NN}_k^j \mu_{k-1}^j}{\sum_{i=1}^{|{\mathbb M}|}{\NN}_k^i \mu_{k-1}^i},\delta\}$;
\ENDFOR
\FOR{$j \in {\mathbb M}$}
\STATE $\mu_k^j = \frac{\bar{\mu}_k^j}{\sum_{i=1}^{|{\mathbb M}|}\bar{\mu}_k^i}$;
\ENDFOR
\STATE Set $\hat{j}_k = \argmax_j \mu_k^j $;
%\\ $\triangleright$ \textbf{\textit{Significance of attacks}}
\STATE Obtain $\chi^2_{|\hat{d}_{1,k}^{\hat{j}_k}|}(\alpha_1)$ and $\chi^2_{|\hat{d}_{2,k-1}^{\hat{j}_k}|}(\alpha_2)$;
%(If $|\hat{d}_{1,k}^{\hat{j}(t)}|=1$ or $|\hat{d}_{2,k-1}^{\hat{j}(t)}|=1$, then $\chi^2_{p=|\hat{d}_{1,k}^{\hat{j}(t)}|}(\alpha_1)=z_{\alpha_1}^2$, $\chi^2_{p=|\hat{d}_{2,k-1}^{\hat{j}(t)}|}(\alpha_2)=z_{\alpha_2}^2$, respectively);
\IF{($\hat{d}_{1,k}^{\hat{j}_k})^T (P_k^{d_1,\hat{j}_k})^{-1}\hat{d}_{1,k}^{\hat{j}_k}<\chi^2_{|\hat{d}_{1,k}^{\hat{j}_k}|}(\alpha_1)$ \textbf{and} $(\hat{d}_{2,k-1}^{\hat{j}_k})^T (P_{k-1}^{d_2,\hat{j}_k})^{-1}\hat{d}_{2,k-1}^{\hat{j}_k}<\chi^2_{|\hat{d}_{2,k-1}^{\hat{j}_k}|}(\alpha_2)$}
\STATE Set $\hat{j}_k^{true}$ as signal attack-free mode of $\hat{j}_k$;
\STATE $\hat{d}_{1,k}^{\hat{j}_k}=\hat{d}_{2,k-1}^{\hat{j}_k}=0$;
\ENDIF 
\STATE \textbf{Return:} 
%$ \hat{j}(t)=\hat{j}_k^{true} {\rm \ for \ } t \in (t_{k-1},t_{k}], \hat{x}(t_k)=\hat{x}_{k|k}^{\hat{j}_k}, t \in  (t_{k-1},t_k], \hat{d}_1(t)=\hat{d}_{1,k}^{\hat{j}_k} {\rm \ for \ } t \in [t_{k},t_{k+1}), \hat{d}_2(t)=\hat{d}_{2,k-1}^{\hat{j}_k} {\rm \ for \ } t \in (t_{k-1},t_k].$ 
\begin{align*}
&\hat{j}(t)=\hat{j}_k^{true} {\rm \ for \ } t \in (t_{k-1},t_{k}],\nnum\\
&\hat{x}(t)=\hat{x}_{k|k}^{\hat{j}_k}, t \in (t_{k-1},t_k],\nnum\\
&\hat{d}_1(t)=\hat{d}_{1,k}^{\hat{j}_k} {\rm \ for \ } t \in [t_{k},t_{k+1}),\nnum\\
&\hat{d}_2(t)=\hat{d}_{2,k-1}^{\hat{j}_k} {\rm \ for \ } t \in (t_{k-1},t_k].
\end{align*}
\ENDFOR
\end{algorithmic}
\end{algorithm}

Now consider the NISME (Algorithm~\ref{algo_mode}) which is derived in Section~\ref{sec:MM}.
The NISME runs the NISE for each mode $j \in {\mathbb M}$ in parallel to generate the state and attack vector estimates along with the priori probability for each mode (line 4).
After then, the algorithm identifies the most likely mode (lines 6-11). 
By the Bayes' theorem, the posteriori probability $\mu_k^j$ is updated by a linear combination of the priori probabilities (line 7).
It is not desirable that some mode probabilities vanish over time because the true modes might be time-varying.
A lower bound $\delta$ is adopted in line 7 to prevent the vanishment of the mode probabilities.
After the lower bound is applied, the mode probability is normalized in line 10.
The mode with the largest posteriori probability $\mu_k^j$ is chosen as a current mode (line 12), and
the attack vectors of the current mode are tested by Chi-square hypothesis tests (p.354 in~\cite{papoulis2002probability}) with significance levels $\alpha_1, \alpha_2$ to determine whether they are statistically significant or not (line 14). Specifically, we have the following null-hypothesis and alternative hypothesis
\begin{align*}
&{\mathcal H}_0: d_{1,k} = 0 {\rm \ and \ } d_{2,k-1}=0\nnum\\
&{\mathcal H}_1: d_{1,k} \neq 0 {\rm \ or \ } d_{2,k-1} \neq 0
\end{align*}
with samples $\hat{d}_{1,k}^{\hat{j}_k}$ and $\hat{d}_{2,k-1}^{\hat{j}_k}$.
Chi-square value is presented as $\chi_{df}^2(\alpha)$ where $df$ and $\alpha$ are the degree of freedom and significance level, respectively.
If it is not statistically significant, the algorithm chooses the signal attack-free mode as a current mode.
The corresponding state and attack vector estimates are returned (line 18).
Due to limited measurements over the continuous-time dynamic system model,
we use the approximation that the attack vector estimates are constants during a sampling period, in lines 2,4,10 of the NISE, and lines 18 of the NISME. We, however, will consider approximation errors in the analysis.

%The mode with the largest a posteriori $\mu_k^j$ is chosen as a current mode and the corresponding state and attack vector estimates are returned (line 13).
%It should be emphasized that the attack vector estimate is approximated as a constant (i.e., piece-wise constant) during a sampling period in lines 2,4,10 of the NISE and lines 12,18 of the NISME. This is because our measurements over the continuous dynamic system are limited to samling instants.

\section{Analysis}\label{sec:analysis0}
In this section, we analyze that the state and attack vector estimation errors of the proposed estimator satisfy PESp-like properties.
The properties in this section are guaranteed for the true mode when the hidden mode is fixed but remains unknown.

Recall that we use the approximation that the attack vector estimates are constants during a sampling period.
In addition, we assume that the attack vectors are continuous in a sampling period, and their gradients are uniformly bounded.
\begin{assumption}
Attack vector $d(t)$ is continuous, and its slope is uniformly bounded; i.e., there exists a positive constant $\bar{d}$ such that
\begin{align*}
\sup_{t_1,t_2 \geq 0}\|(d(t_1)-d(t_2))/(t_1-t_2)\| \leq \bar{d}.
\end{align*}
\label{CD_asm1}
\end{assumption}
Some required notations are introduced below.
It is obvious that $w_k$ shown in~\eqref{DD001} is uncorrelated with $v_{1,k}$, $v_{2,k}$, or $v_{3,k}$.
Noise vectors $v_{1,k}$, $v_{2,k}$, and $v_{3,k}$ are also uncorrelated with each others because
\begin{align*}
R_{1,k}&={\mathbb E}[v_{1,k}v_{1,k}^T]= T_{1,k}R_kT_{1,k}^T> 0\nnum\\
R_{2,k}&={\mathbb E}[v_{2,k}v_{2,k}^T]= \bar{T}_{2,k}T_{2}R_kT_{2}^T\bar{T}_{2,k}^T> 0\nnum\\
R_{3,k}&={\mathbb E}[v_{3,k}v_{3,k}^T]= \bar{T}_{1,k}T_{2}R_kT_{2}^T\bar{T}_{1,k}^T> 0\nnum\\
R_{12,k}&={\mathbb E}[v_{1,k}v_{2,k}^T]= (U_{1}^TR_kU_{2}^T-U_{1}^TR_kU_{2}^T\nnum\\
&\times (U_{2}^TR_kU_{2}^T)^{-1}U_{2}^TR_kU_{2}^T)\bar{T}_{2,k}^T=0\nnum\\
R_{13,k}&={\mathbb E}[v_{1,k}v_{3,k}^T]= (U_{1}^TR_kU_{2}^T-U_{1}^TR_kU_{2}^T\nnum\\
&\times (U_{2}^TR_kU_{2}^T)^{-1}U_{2}^TR_kU_{2}^T)\bar{T}_{1,k}^T=0
\end{align*}\begin{align*}
R_{23,k}&={\mathbb E}[v_{2,k}v_{3,k}^T]= \bar{U}_{1,k}^TT_{2}R_kT_{2}^T\bar{U}_{2,k}^T-\bar{U}_{1,k}^TR_k\bar{U}_{2,k}^T\nnum\\
&\times (\bar{U}_{2,k}^TR_k\bar{U}_{2,k}^T)^{-1}\bar{U}_{2,k}^TT_{2}R_kT_{2}^T\bar{U}_{2,k}^T=0.
\end{align*}
%For the same reason, it holds that ${\mathbb E}[C_{2,k}w_{k-1}x_k^TC_{3,k}^T]=0$.
We introduce the linearization errors $\phi_k$, $\psi_{1,k}$, $\psi_{2,k}$, and $\psi_{3,k}$ for the stability analysis:
\begin{align}
& f(x_{k},u_k,d_{1,k},w_k',t_k)- f(\hat{x}_{k|k},u_{k},\hat{d}_{1,k},0,t_{k})\nnum\\
&= A_k\tilde{x}_{k|k}+ G_{1,k}\tilde{d}_{1,k}+w_k+ \phi_k(\hat{x}_{k|k},x_k,u_k,w_k',v_k')\nnum\\
&h_1(x_k,u_k,v_k',t_k)-h_1(\hat{x}_{k|k},u_k,0,t_k)\nnum\\
&=C_{1,k}\tilde{x}_{k|k}+v_{1,k}+\psi_{1,k}(\hat{x}_{k|k},x_k,u_k,v_k')\nnum\\
&h_2(u_k,v_k',t_k)-h_2(u_k,0,t_k)\nnum\\
&=v_{2,k}+\psi_{2,k}(u_k,v_k')\nnum\\
&h_3(x_k,u_k,v_k',t_k)-h_3(\hat{x}_{k|k-1},u_k,0,t_k)\nnum\\
&=C_{3,k}\tilde{x}_{k|k-1}+v_{3,k}+\psi_{3,k}(\hat{x}_{k|k-1},x_k,u_k,v_k')
\label{Non_er}
\end{align}
where $\phi_k$ is a function of $\tilde{d}_{1,k}$ and $\tilde{d}_{1,k}$ is a function of $\tilde{x}_{k|k}$ and $v'_k$. 
We omit the arguments of the linearization errors in the rest of the paper.
The linearization errors can not be used to estimate the states and attack vectors, since they include unknown variables $x_k$, $w'_k$, and $v'_k$.

%\subsection{Stability}
Now consider the following assumptions.
\begin{assumption}
There exist positive constants $\bar{a}'$, $\bar{c}_3$, $\underline{q}'$, and $\underline{r}_{3}$ such that the following holds for $\forall k$:
\begin{align*}
&\|A_k\| \leq \bar{a}', \ \|C_{3,k}\| \leq \bar{c}_3, \ \underline{q}' \leq Q_k, \ \underline{r}_{3}I \leq R_{3,k}.
\end{align*}
If $Rank(\Sigma) \neq 0$, there exist positive constants
$\bar{c}_{1}$, $\bar{g}_1$, and $\bar{m}_1$ such that the following holds for $\forall k$:
 \begin{align*}
\|C_{1,k}\| \leq \bar{c}_{1}, \ \|G_{1,k}\| \leq \bar{g}_1, \ \|\Sigma^{-1}\| \leq \bar{m}_1.
\end{align*}
If $Rank(\bar{\Sigma}_k) \neq 0$, there exist positive constants $\bar{c}_{2}$, $\underline{g}_2$, $\bar{g}_2$, $\underline{m}_2$, $\bar{m}_2$, and $\underline{r}_2$ such that the following holds for $\forall k$:
\begin{align*}
&\|C_{2,k}\| \leq \bar{c}_{2}, \underline{g}_2 \leq \|G_{2,k}\| \leq \bar{g}_2, \underline{m}_2 \leq \|\bar{\Sigma}_k^{-1}\| \leq \bar{m}_2, \nnum\\
& \underline{r}_2 I \leq R_{2,k}.
\end{align*}
\label{asm_1}
\end{assumption}
\begin{assumption}
For any $\epsilon_{\phi}, \epsilon_{\psi_{1}}, \epsilon_{\psi_2}, \epsilon_{\psi_3}>0$, there exists $\delta>0$ such that 
\begin{align*}
&\|\phi_k\| \leq \epsilon_{\phi}\|x_k-\hat{x}_{k|k}\|^2, \ \ \|\psi_{1,k}\| \leq \epsilon_{\psi_{1}}\|x_k -\hat{x}_{k|k}\|^2\nnum\\
&\|\psi_{2,k}\| \leq \epsilon_{\psi_2}\|x_k-\hat{x}_{k|k}\|^2, \ \ \|\psi_{3,k}\| \leq \epsilon_{\psi_3}\|x_k-\hat{x}_{k|k}\|^2
\end{align*}
hold for all $\|x_k-\hat{x}_{k|k}\| \leq \delta$ and $k\geq0$. For any $\epsilon_{\rho}>0$, there exists $\delta_{\rho}>0$ such that
$\|\rho_{k}\| \leq \epsilon^2 \epsilon_{\rho}$ for all $\epsilon \leq \delta_{\rho}$ and $k\geq 0$, where $\rho_k$ is defined in~\eqref{E001}.
\label{asm_2}
%\footnote{$\phi_k$ is a function of $\epsilon$.}
%\item For all $\epsilon_{\phi}, \epsilon_{\psi}>0$, there exist $\delta_{\phi}, \delta_{\psi}>0$ such that\begin{align*} \|\phi_k\| \leq \epsilon_{\phi}\|x_k-\hat{x}_{k|k}\|^2, \ \ \|\psi_k\| \leq \epsilon_{\psi}\|x_k-\hat{x}_{k|k-1}\|^2 \end{align*} for all $\|x_k-\hat{x}_{k|k}\| \leq \delta_{\phi}$ and $\|x_k-\hat{x}_{k|k-1}\| \leq \delta_{\psi}$.
\end{assumption}

\begin{assumption}
There exist positive constants $\underline{p}$ and $\bar{p}$ such that $\underline{p}I \leq P_{k}^x \leq \bar{p}I$ for $\forall k$. 
\label{asm_11}
\end{assumption}
Under the assumptions, we can guarantee PESp-like properties for the estimation errors of states and attack vectors.
\begin{theorem}
Consider the NISE, provided that Assumptions~\ref{CD_asm_m}, \ref{CD_asm1}, \ref{asm_1}, \ref{asm_2}, and~\ref{asm_11} hold.
For any $\gamma \in (0,1)$, there exist positive constants $\alpha_x$, $\alpha_{d_1}$, $\alpha_{d_2}$, $b_x$, $b_{d_1}$, $b_{d_2}$, $c_x$, $c_{d_1}$, $c_{d_2}$, $\underline{\delta}$, $\bar{q}'$, $\bar{r}_1$, $\bar{r}_2$, $\bar{r}_3$, and $\bar{\epsilon}$ such that, 
if $Q_k \leq \bar{q}' I$, $R_{1,k} \leq \bar{r}_1 I$, $R_{2,k} \leq \bar{r}_2 I$, $R_{3,k} \leq \bar{r}_3 I$, and $\epsilon \leq \bar{\epsilon}$, then the following properties hold:
\begin{align*}
P(\|\tilde{x}_{k|k}\| &< \alpha_x e^{-b_x k }\|\tilde{x}_{0|0}\| + c_x) \geq 1-\gamma,\nnum\\
P(\|\tilde{d}_1(t)\| &< \alpha_{d_1} e^{-b_{d_1} t }\|\tilde{x}_{0|0}\| + c_{d_1}) \geq 1-\gamma,\nnum\\
P(\|\tilde{d}_2(t)\| &< \alpha_{d_2} e^{-b_{d_2} t }\|\tilde{x}_{0|0}\| + c_{d_2}) \geq 1-\gamma\nnum\\
\end{align*}
for all $\|\tilde{x}_{0|0}\| \leq \underline{\delta}$, $k\geq 0$ and $t\geq0$.
\label{CD_stability_the1}
\end{theorem}
Theorem~\ref{CD_stability_the1} is formally proven in Section~\ref{AP_CD_Th_proof}.
A sufficient condition for Assumptions~\ref{asm_1} is that functions $f$ and $h$ satisfy bi-Lipscitz continuity (p.10 in~\cite{ostrovskii2013metric}), and covariance matrices of noise vectors $w_k'$ and $v_k'$ are uniformly lower bounded by positive definite matrices.
Assumptions~\ref{asm_2} can be verified by checking if functions $f$ and $h$ in system~\eqref{CD001a} satisfy Holder continuity (p.136 in~\cite{truesdell2004non}) with exponent $2$.
%Assumptions~\ref{asm_1} and~\ref{asm_2} can be verified off-line by checking if functions $f$ and $h$ satisfy bi-Lipscitz continuity (sufficiency for Assumptions~\ref{asm_1}) and holder continuity with exponent $2$ (sufficiency for Assumptions~\ref{asm_2}).}
One sufficient condition for Assumption~\ref{asm_11} is the uniform observability of the pair $(C_{3,k},\bar{A}_k)$, as shown in Lemma~\ref{Bounded_2}. We recall that $C_{3,k}$ and $\bar{A}_k$ are defined in~\eqref{CD002} and \eqref{CD004.7}.
\begin{lemma}
Consider the NISE.
Under Assumptions~\ref{CD_asm_m},~\ref{asm_1}, and~\ref{asm_2}, if the pair $(C_{3,k},\bar{A}_k)$ is uniformly observable and $P_{0}^x \geq 0$, then there exist positive constants $\underline{p}$ and $\bar{p}$ such that $\underline{p}I \leq P_k^x \leq \bar{p}I$ for all $k \geq 0$.
\label{Bounded_2}
\end{lemma}
Lemma~\ref{Bounded_2} is formally proven in Section~\ref{AP_CD_Th_proof2}.
%Lemma~\ref{Bounded_2} does not depend on the additional errors induced by the piece-wise constant attack vector approximation because the covariance matrices are updated only at sampling instants.
Lemma~\ref{Bounded_2} does not require Assumption~\ref{CD_asm1} because the update rule of error covariance matrices depend only on known system matrices and covariance matrices of noises.
\begin{remark}
Corresponding parameters for PESp-like properties in Theorem~\ref{CD_stability_the1} can be obtained in the proof.
For ease of presentation, we omit the procedure to find these parameters. 
\oprocend
\end{remark}
%\begin{remark}Theorem~\ref{CD_stability_the1} and Lemma~\ref{Bounded_2} as well as Assumptions~\ref{asm_1},~\ref{asm_2} and~\ref{asm_11} are extensions of the discrete time extended Kalman filter~\cite{kluge2010stochastic} by integrating the attack vectors and output discretization. When there are no attack vectors with zero discretization errors, they reduce to Theorem 3.1, 4.5, and Assumption 3.1 and 4.3 in~\cite{kluge2010stochastic} of the extended Kalman filer. To our best knowledge, this is the first time to analyze the properties of the NISE. \oprocend \end{remark}

\section{Discussion: Reducing mode number}\label{sec:DiscR}
When there are ${s+m}$ signal attack locations, mode set ${\mathbb M}^A$ includes all the combinations of the attack locations; i.e., $|\mathbb M^A| = 2^{s+m}$.
As the number of signal attack location increases, computational complexity increases exponentially.
We, in this section, discuss how to alleviate computational complexity by reducing the number of modes induced by ${\mathbb M}^A$, and how to estimate the true mode from the estimation results of a reduced mode set for hidden-mode switched linear systems.
Finding a reduced mode set presented in Section~\ref{sec:MR}, and finding a true mode presented in Section~\ref{sec:truemode} are both on-line procedures. The former is conducted before running the NISE, and the latter will replace the hypothesis test (line 13-17) in the NISME. 
The complete algorithm for the NISME with reduced mode set is presented in Section~\ref{sec:truemode}, where we consider the special case with $|{\mathbb M}^I|=1$.
In Section~\ref{sec:NISE_re}, we extend the results to any ${\mathbb M}^I$.

Let us consider the special case where the dynamic for each mode $j$ is linear and switching, and $|{\mathbb M}^I|=1$:
\begin{align}
\dot{x}(t) &= A(t) x(t) +B(t)(u(t) + S^{j(t)} d(t))+w(t),\nnum\\
& \ \ \quad \quad \quad \quad \quad \quad \quad \quad \quad  \quad \quad \quad \quad \quad \quad \quad x(t) \in {\CC}^{j(t)}\nnum\\
(x(t),&j(t))^+=\Omega(x(t),j(t)), \quad \quad \quad \quad \quad \quad \quad x(t) \in {\DD}^{j(t)}\nnum\\
y_k &= C_k x_k +D_k u_k +H^{j_k} d_{k}+v_k
\label{e_dis}
\end{align}
where %$G^{j(t)}(t)=B(t)L^{j(t)}$ and
$S^{j(t)}$, and $H^{j_k}$ are defined in~\eqref{CD001a}.
We remind $K^j = [(S^j)^T,(H^j)^T]^T$, and define ${\mathbb K}^{j}_1 \triangleq \{i =1,\cdots,n+m|K^j(i,i) =1\}$.

%The number of possible modes is $2^{m}$, but their performance might be poor.

\subsection{Mode reduction}\label{sec:MR}
Mode reduction is based on the following two ideas. Firstly, we maintain the modes such that uniform observability over finite time-horizon $[c_1,c_2]$ (Definition~\ref{defi:uniform2}) holds. Secondly, we remove modes whose assumptions on attack locations are strictly restrictive than those of others.
\begin{definition}
The pair $(C_k,A_k)$ is uniformly observable over $[c_1,c_2]$ if and only if there exist positive constants $a,b$, and $l<c_2-c_1$ such that
$a I \leq {\mathcal M}_{k+l,k} \leq b I$ for $k = c_1,c_1+1,\cdots, c_2-l$.% where $M_{k+l,k} \triangleq \sum_{i = k}^l \Phi_{k+1,i+1}C_kC_k^T \Phi_{k+1,i+1}^T$ is the observability gramian and $\Phi_{k,l}= U_kU^{-1}_k$ is the state transition matrix with fundamental matrix $U_{k+1}=A_kU_k$.
\label{defi:uniform2}
\end{definition}

The first idea is motivated by the sufficient condition for Lemma~\ref{Bounded_2}; uniform observability.
To check uniform observability, the defender is required to have information on pairs $(C_{3,k}^j,\bar{A}_k^j)$ for $k =0,1,\cdots$. 
This information is hard to gather at initial time. We instead adopt an approximation, Definition~\ref{defi:uniform2}, which only requires the system matrices for a few next steps. Uniform observability over $[c_1,c_2]$ reduces to uniform observability as $c_2 \rightarrow \infty$ with fixed $c_1$.

To justify the second idea, consider a pair of modes $j, j' \in {\mathbb M}^A$ such that ${\mathbb K}^{j'}_1 \subset {\mathbb K}^{j}_1$, and $(C_{3,k}^j,\bar{A}_k^j)$, $(C_{3,k}^{j'},\bar{A}_k^{j'})$ are uniformly observable over $[c_1,c_2]$.
The relation ${\mathbb K}^{j'}_1 \subset {\mathbb K}^{j}_1$ indicates that mode $j'$ imposes a more restrictive assumption on attack locations than mode $j$.
In this sense, mode $j'$ is said to be redundant and it could be ruled out to reduce computational complexity.

Intuitively speaking, the above ideas allow the minimal number of modes to provide the same attack capability as the power set.
The reduced mode set is defined by ${\mathbb M}^{d}_{[c_1, c_2]} = \{j \in {\mathbb M}^{ob}_{[c_1, c_2]} |\nexists j' \in {\mathbb M}^{ob}_{[c_1, c_2]} {\rm \ s.t. \ } {\mathbb K}^j_1 \subset {\mathbb K}^{j'}_1\}$ where
${\mathbb M}^{ob}_{[c_1, c_2]} \triangleq \{j \in {\mathbb M}^A | (C_{3,k}^{j},\bar{A}^j_k)$ is uniformly observable over $[c_1,c_2]\}$.
It can be found by Algorithm~\ref{algo3} where ${\mathbb M}_{i}^A \triangleq \{j \in {\mathbb M}^A | |{\mathbb K}_1^j|=i\}$.
%is upper bounded by $C_{s+m}^{\left \lfloor{\frac{s+m}{2}}\right \rfloor}$. This is summarized in Lemma~\ref{lemm002} and Section~\ref{ap_disc} presents its proof.
Without the mode reduction, the worst upper bound of $|{\mathbb M}^A|$ is $2^{s+m}$, but ${\mathbb M}^{d}_{[c_1, c_2]}$ could be as low as $1$ (see case study 2).
It is worthy to emphasize that the defender needs to know $(C_{3,k}^j,\bar{A}_k^j)$ over $[c_1,c_2]$ in order to verify uniform observability in the interval.

\begin{algorithm} \caption{Mode reduction (finding ${\mathbb M}^{d}_{[c_1, c_2]}$).} \label{algo3}
\begin{algorithmic}[1]
\REQUIRE ${\mathbb M}^A$, $A_k$, $C_k$ for $k =c_1,c_1+1,\cdots,c_2$ (or corresponding
 $A(t)$, $C(t)$), $G^j$, $H^j$ for $\forall j \in {\mathbb M}^A$;
\STATE ${\mathbb M}^{d}_{[c_1, c_2]}=\emptyset$;
\FOR{$i=s+m:1$}
\FOR{$j \in {\mathbb M}_{i}^A$}
\IF{$(C_{3,k}^j,\bar{A}_k^j)$ is uniformly observable over $[c_1,c_2]$, \textbf{and}
$\nexists j' \in {\mathbb M}^d_{[c_1, c_2]} {\rm \ s.t. \ } {\mathbb K}^j_1 \subset {\mathbb K}^{j'}_1$}
\STATE ${\mathbb M}^{d}_{[c_1, c_2]}={\mathbb M}^{d}_{[c_1, c_2]} \cup \{j\}$;
\ENDIF
\ENDFOR
\ENDFOR
\ENSURE ${\mathbb M}^{d}_{[c_1, c_2]}$.
\end{algorithmic}
\end{algorithm}

%\begin{lemma} It holds that $|{\mathbb M}^{d}_{[c_1, c_2]}| \leq C_{s+m}^{\left \lfloor{\frac{s+m}{2}}\right \rfloor}$. \label{lemm002} \end{lemma}

\begin{remark}
Mode reduction Algorithm~\ref{algo3} is not applicable to nonlinear systems because uniform observability is determined by where linearization is performed. This information cannot be obtained in advance.
\oprocend
\end{remark}

\subsection{True mode estimation}\label{sec:truemode}
We discuss how to estimate the true mode from the outputs of the NISME under the reduced mode set.
It might be noticed that the reduced mode set might not include true modes, since some of the modes are removed.
Given mode estimate $\hat{j}(t)$ from the reduced mode set ${\mathbb M}^d_{[c_1, c_2]}$, the idea is to conduct two-tailed $z$-test~\cite{papoulis2002probability} for each attack location $i \in {\mathbb K}_1^{\hat{j}(t)}$ to determine whether the attack size is statistically significant. To be specific, we test the null hypothesis that $i^{th}$ elements of $d_{1,k}$ or $d_{2,k-1}$ are zero:
\begin{align*}
{\mathcal H}_0: d_{1,k}(i-s) = 0 {\rm \ if \ } i > s, \ {\mathcal H}_0: d_{2,k-1}(i) = 0 {\rm \ if \ } i \leq s
\end{align*}
and $\hat{d}_{1,k}^{\hat{j}_k}(i-s)$ or $\hat{d}_{2,k-1}^{\hat{j}_k}(i)$ are regarded as samples.
If the null hypothesis is rejected, then we accept alternative hypothesis
\begin{align*}
{\mathcal H}_1: d_{1,k}(i-s) \neq 0 {\rm \ if \ } i > s, \  {\mathcal H}_1: d_{2,k-1}(i) \neq 0 {\rm \ if \ } i \leq s
\end{align*}
i.e., there exists an attack on $i^{th}$ location.
Algorithm~\ref{algo4} presents the pseudo code for true mode estimation.
$z$-value is presented as $z(\alpha)$ where $\alpha$ is the significance level.
Hypothesis tests are conducted in lines 6 and 13 for actuator attacks, and sensor attacks, respectively.
%We could conduct hypothesis tests for attack vectors element-wise because $K_g$ and $K_h$ are diagonal matrices.

\subsection{NISME with reduced mode set}\label{sec:NISE_re}
Algorithm~\ref{algo5} shows the NISME with reduced mode sets.
The core of the algorithm is identical to that of the NISME, and some differences are explained as follows.
We apply the mode reduction technique and true mode estimation technique to each $i \in {\mathbb M}^I$.
The algorithm calculates reduced mode set every $\TT$ steps for every $i \in {\mathbb M}^I$ (line 4).
This requires the defender to have knowledge on system matices for next $\TT-1$ steps.
Based on the fact that the reduced mode set might not include the true mode, we test attack vectors element-wise to identify the true mode (line 19).
As ${\TT}$ decreases in Algorithm~\ref{algo5}, lesser knowledge on system matrices is required, but computational complexity induced by Algorithm~\ref{algo3} increases. When system~\eqref{e_dis} is time-invariant, ${\TT}=\infty$.
%If a reduced mode set is used, Algorithm~\ref{algo4} will replace line 14-19 in Algorithm~\ref{algo_mode}, and the mode outputs of the NISME be $\hat{j}^{true}(t)$.

\begin{algorithm} \caption{Mode estimation with mode reduction} \label{algo4}
\begin{algorithmic}[1]
\REQUIRE $\hat{j}_k$, $\hat{d}_{1,k}^{\hat{j}_k}$, $\hat{d}_{2,k-1}^{\hat{j}_k}$, $P_k^{d_1,\hat{j}_k}$, $P_{k-1}^{d_2,\hat{j}_k}$, $\alpha_1$, $\alpha_2$ (significance levels);
\STATE Obtain $z$-values $z(\alpha_1)$ and $z(\alpha_2)$ from $z$-test table;
\STATE %$\hat{d}_1=\hat{d}_1(t)$, $\hat{d}_2=\hat{d}_2(t)$
$K^{\hat{j}^{true}_k}=0^{(s+m) \times (s+m)}$;
\STATE $l_1=l_2=1$;
\FOR{$i \in {\mathbb K}^{\hat{j}_k}_1$}
\IF{$i \leq s$}
\IF{$\frac{|\hat{d}_{2,k-1}^{\hat{j}_k}(l_2)|}{\sqrt{P^{d_2,\hat{j}_k}_{k-1}(l_2,l_2)}}> z(\alpha_2)$}
\STATE $K^{\hat{j}^{true}_k}(i,i)=1$, $\hat{d}_{2,k-1}^{\hat{j}_k^{true}}(l_2)=\hat{d}_{2,k-1}^{\hat{j}_k}(l_2)$;
\ELSE
\STATE $\hat{d}_{2,k-1}^{\hat{j}_k^{ture}}(l_2)=0$;
\ENDIF
\STATE $l_2=l_2+1$;
\ELSE
\IF{$\frac{|\hat{d}_{1,k}^{\hat{j}_k}(l_1)|}{\sqrt{P_k^{d_1,\hat{j}_k}(l_1,l_1}}> z(\alpha_1)$}
\STATE $K^{\hat{j}^{true}_k}(i,i)=1$, $\hat{d}_{1,k}^{\hat{j}^{true}_k}(l_1)=\hat{d}_{1,k}^{\hat{j}_k}(l_1)$;
\ELSE
\STATE $\hat{d}_{1,k}^{\hat{j}_k^{true}}(l_1)=0$;
\ENDIF
\STATE $l_1=l_1+1$;
\ENDIF
\ENDFOR
\STATE Obtain $\hat{j}^{true}_k$ for corresponding $K^{\hat{j}^{true}_k}$;
\ENSURE $\hat{j}^{true}_k$, $\hat{d}_{1,k}^{\hat{j}^{ture}_k}$, $\hat{d}_{2,k-1}^{\hat{j}^{ture}_k}$.
\end{algorithmic}
\end{algorithm}

\begin{algorithm} \caption{NISME with reduced mode set} \label{algo5}
\begin{algorithmic}[1]
\REQUIRE $\hat{x}_{0|0}^j={\mathbb E}[x_0]$, $P_0^{x,j}=P_0$, $\mu_0^j=\frac{1}{|{\mathbb M}|}$,
$\hat{d}_{1,0}^j=(\Sigma^j)^{-1}(z_{1,0}^j-T_{1,0}^jh(\hat{x}_{0|0}^j,u_0,0,j,0))$ for $\forall j \in {\mathbb M}$; Choose $0<\delta \ll \frac{1}{|{\mathbb M}|}$, $0<\alpha_1<1$, $0<\alpha_2<1$ (significance levels), and ${\TT}$;
\FOR{$q=0:N$}
\STATE $c_1=1+q\TT$, $c_2=(q+1)\TT$, ${\mathbb M}^D_{[c_1, c_2]} = \emptyset$;
\FOR{$i \in {\mathbb M}^I$}
\STATE Run Algorithm~\ref{algo3} with (${\mathbb M}^{A}$, $A_k^i$, $C_k^i$ for $k =c_1,c_1+1,\cdots,c_2$, $G^j$, $H^j$ for $\forall j \in {\mathbb M}^A$) to obtain ${\mathbb M}^{d,i}_{[c_1, c_2]}$;
\STATE ${\mathbb M}^D_{[c_1, c_2]} = {\mathbb M}^D_{[c_1, c_2]} \cup {\mathbb M}^{d,i}_{[c_1, c_2]}$;
\ENDFOR
\FOR{$k=c_1:c_2$}
\STATE Read sensor output $y_k$, and control input $u(t)$ for $t \in [t_{k-1},t_{k}]$;
\FOR{$j \in {\mathbb M}^D_{[c_1, c_2]}$}
\STATE Run the NISE with input $(j,\hat{x}_{k-1|k-1}^j,\hat{d}_{1,k-1}^j,$ $P_{k-1}^{x,j}, P_{k-1}^{d_1,j}, P_{k-1}^{xd_1,j},y_k,u(t) {\rm \ for \ } t \in [t_{k-1},t_k])$ to generate output ($\hat{x}_{k|k}^j$, $\hat{d}_{1,k}^j$, $\hat{d}_{2,k-1}^j$, $P_{k}^{x,j}$, $P_{k-1}^{d_2,j}$, $P_{k}^{d_1,j}$, $P_{k}^{xd_1,j}$, ${\NN}_k^j$);
\ENDFOR \\ 
$\triangleright$ \textbf{\textit{Mode estimator}}
\FOR{$j \in {\mathbb M}^D_{[c_1, c_2]}$}
\STATE $\bar{\mu}_k^j = \max\{\frac{{\NN}_k^j \mu_{k-1}^j}{\sum_{i=1}^{|{\mathbb M}^D_{[c_1, c_2]}|}{\NN}_k^i \mu_{k-1}^i},\delta\}$;
\ENDFOR
\FOR{$j \in {\mathbb M}^D_{[c_1, c_2]}$}
\STATE $\mu_k^j = \frac{\bar{\mu}_k^j}{\sum_{i=1}^{|{\mathbb M}^D_{[c_1, c_2]}|}\bar{\mu}_k^i}$;
\ENDFOR
\STATE Set $\hat{j}_k = \argmax_j \mu_k^j $;
\STATE Run Algorithm~\ref{algo4} with ($\hat{j}_k$, $\hat{d}_{1,k}^{\hat{j}_k}$, $\hat{d}_{2,k-1}^{\hat{j}_k}$, $P_k^{d_1,\hat{j}_k}$, $P_{k-1}^{d_2,\hat{j}_k}$, $\alpha_1$, $\alpha_2$) to obtain ($\hat{j}^{true}_k$, $\hat{d}_{1,k}^{\hat{j}_k^{ture}}$, $\hat{d}_{2,k-1}^{\hat{j}_k^{ture}}$);
\STATE \textbf{Return:} 
\begin{align*}
&\hat{j}(t) = \hat{j}^{true}_k {\rm \ for \ } t \in (t_{k-1},t_{k}],\nnum\\
&\hat{x}(t)=\hat{x}_{k|k}^{\hat{j}_k}, t \in (t_{k-1},t_k],\nnum\\
&\hat{d}_1(t)=\hat{d}_{1,k}^{\hat{j}_k^{ture}} {\rm \ for \ } t \in [t_{k},t_{k+1}),\nnum\\ &\hat{d}_2(t)=\hat{d}_{2,k-1}^{\hat{j}_k^{ture}} {\rm \ for \ } t \in (t_{k-1},t_k].
\end{align*}
\ENDFOR
\ENDFOR
\end{algorithmic}
\end{algorithm}
\section{Numerical simulation}\label{sec:nu_sim}
\begin{figure}[h]
  \centering
  \includegraphics[width = \linewidth]{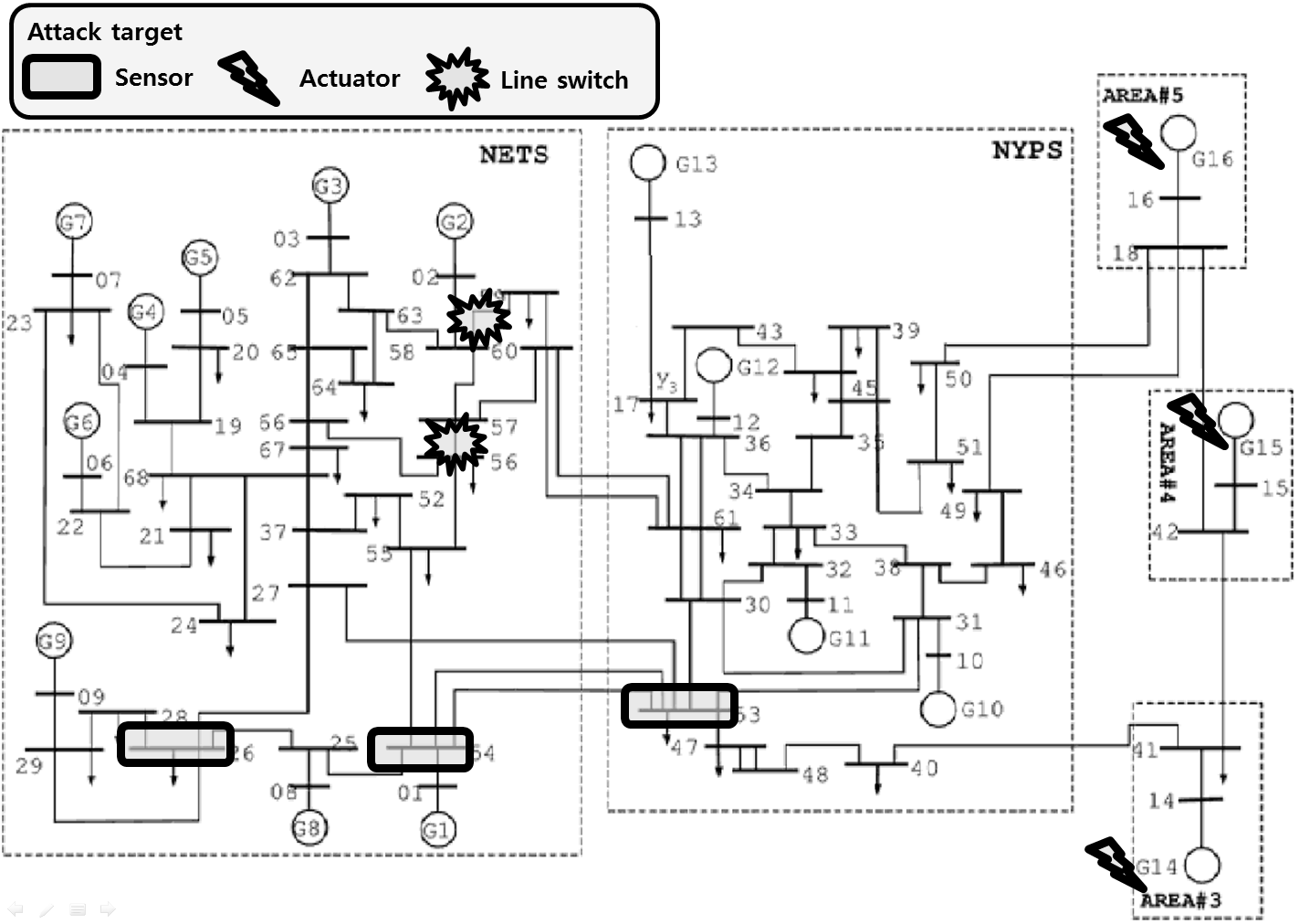}
  \caption{Locations of the attacks (Figure from~\cite{zhang2012flexible}).}\label{TP68}
\end{figure}

In Section~\ref{sec:analysis0}, a set of properties is shown for the true mode under the time-invariant hidden mode. However, we do not have formal guarantees for the cases of time-varying modes.
Moreover, the effectiveness of a set of reduced modes discussed in Section~\ref{sec:DiscR} remains unclear.  
In this section, we will use the IEEE 68-bus test system to empirically illustrate them.
%In this section, we present an illustrative example to show the performance of the proposed NISME.
The NISME is applied to the IEEE 68-bus test system shown in Figure~\ref{TP68}. %~\cite{pal2006robust,rogers2012power}.
In the network, there are $16$ generator buses ($|{\GG}|=16$), and $52$ load buses ($|{\LL}|=52$). Each local bus is described by~\eqref{Po_model} (as~\cite{ilic2010modeling}), and~\eqref{Po_model_output} with $\epsilon=0.01s$. It is assumed that noises $w(t)$ and $v_k$ are zero mean Gaussian with covariance matrices $Q_{i}(t)=0.01^2 I$, and $R_{i,k}=0.01^4 I$.
The parameters are adopted from page 598 in~\cite{Kundar.Balu.Lauby:94}:
$D_i = 1$, and $t_{ij} = 1.5$ for $\forall i \in \VV$.
Angular momentums are $m_i = 10s$ for $i \in \GG$ and a larger value $m_i = 100s$ for load buses $i \in \LL$.
%as suggested in~\cite{zhao2014design}.
Backstepping inspired stabilizing distributed controllers~\cite{HK-MZ:AUTOMATICA15} are applied to the power system.
We choose $\delta=3.3\%$ as a lower bound of probabilities. 
  %Power demand $P_{L_i}(t) = 1+0.1 \sin(0.01 t)$ is known.

The attacker could launch 3 sensor attacks, 3 actuator attacks, and 2 switching attacks described in Figure~\ref{TP68}. %To quantify the size of power line switching attacks, we modify the system model into $\dot{f}_i(t) = -\frac{1}{m_i}\big(D_i f_i(t) + \sum_{l\in {\mathcal S}_i} (P_{il}(t)+d_{c,i,l}(t))+P_{L_i}(t)- (P_{M_i}(t)+d_{a,i}(t))\big)+w_{2,i}(t)$ where power line switching attack vectors $d_{c,i,l}(t)=-d_{c,l,i}(t)$ are introduced in lieu of hidden-mode $j_{il}(t)$.

We consider the attack scenario where the system is under the time-varying attacks: sensor attacks $0.01 \cos(0.12 t)$ for $t = [0,10)$, actuator attacks $0.1-0.6\sin(0.3t)$ for $t = [10,20)$, and switching attacks for $t = [20,30)$.
%The attacker injects sinusoidal sensor attacks $5 \cos(0.12 t)$ and actuator attacks $60 \sin(0.3t)$ and breaks the circuit on the corresponding target buses.
For $t \geq 30$, the system would be attack-free.
%The defender cannot access to the current mode and its transitions.

The goals of case study 1 and 2 are to verify the performance of the NISME for time-varying modes with a regular mode set, and a reduced mode set, respectively.

\textbf{Case study 1:} Consider the following four modes.\\
{Mode 0:} Attack-free.\\
{Mode 1:} Sensors (electrical power outputs) $26, 53$ and $54$ are attacked.\\
{Mode 2:} Actuators $14,15$, and $16$ are attacked.\\
{Mode 3:} Line switches $\{56,57\}$, and $\{59,60\}$ are attacked.

If the sizes of the estimated attack vectors are not statistically significant, mode 0 will be chosen, as described in Algorithm~\ref{algo_mode}.

%either mode 1 or mode 2 can describe this case with zero attack vectors.

\begin{remark}
The power system under the mode 1,2 and 3 satisfies Assumptions~\ref{CD_asm_m},~\ref{CD_asm1},~\ref{asm_1},~\ref{asm_2} and uniform observability condition for all three modes. This is because the system is time invariant and the linearization error is $O(\|\tilde{x}_{k|k}\|^2)$ for $\|\tilde{x}_{k|k}\| \leq 1$.
%If the rest of the assumptions in Theorem~\ref{CD_stability_the1} are satisfied at the initial time, and at the every switching time instants for the corresponding true mode $j(t)$, it is expected that the state estimation errors are PESp for the corresponding true mode NISE.
\oprocend
\end{remark}

%Since the simulated system satisfies all the assumptions, it is expected that, for the real mode, there exists a set of parameters and initial condition such that the state estimation error is exponentially bounded in mean square. 

\begin{figure}[h]
  \centering
  \includegraphics[width = \linewidth]{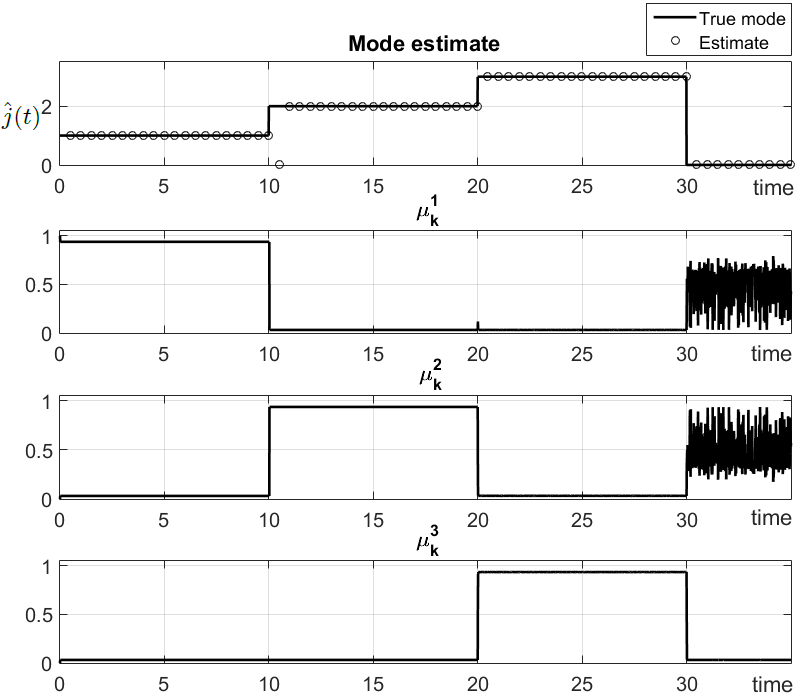}
  \caption{Mode estimates and probabilities of each mode.}\label{Case1_2}
\end{figure}

\begin{figure}[h]
  \centering
  \includegraphics[width = \linewidth]{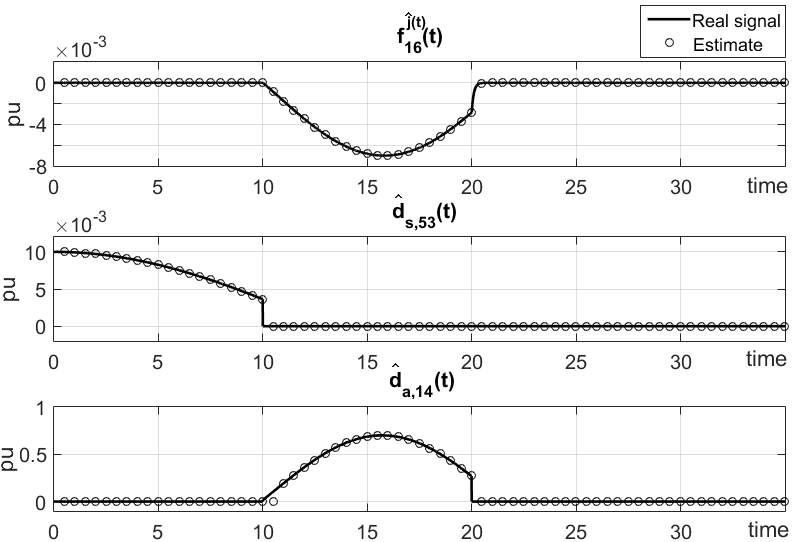}
  \caption{Real signals and estimated signals (top to bottom); (a) state estimation of angular frequency at bus 16; (b) sensor attack estimation of bus 53;  (c) actuator attack estimation on the control input of bus 14.}\label{Case1_1}
\end{figure}

Significance levels $\alpha_1=\alpha_2=0.75$ are applied with corresponding chi-square values $\chi_3^2(0.75)=4.11$.
The mode probabilities and estimation results are shown in Figure~\ref{Case1_2} where the estimates are coincident with the true modes.
Mode estimation is inaccurate near $10$ sec.
This is because the sizes of attack vectors are small (the second and third subfigure in Figure~\ref{Case1_1}) during this time and thus the attack vector estimates are not considered statistically significant.
Mode probabilities of mode 1 and 2 are oscillating for $t>30$, because both modes 1 and 2 with zero attack vectors would represent attack-free mode (mode 0).

The outputs of the NISME and the real attack signals are presented in Figure~\ref{Case1_1}.
The first subfigure indicates that the state estimation errors satisfy PESp-like property. Although the frequency fluctuates due to the actuator attack $t \in [10,20)$, its estimates are accurate.
The additive sinusoidal sensor attack for $t \in [0,10)$ is well estimated as shown in the second subfigure.
The third subfigure shows the estimates and real vectors of actuator attack; the sinusoidal actuator attacks for $t \in [10,20)$ on the control inputs. Around $10$ sec, attack vector estimates are set to zero because $\hat{d}_{a,14}^{j=2}$ is not considered statistically significant.
%The last subfigure presents the estimation result of mode 3; i.e., step signal for $t \in [20,30)$ induced by the power line switching attack.

\textbf{Case study 2:}
There are $6$ potential signal attack locations with $2$ possible switching attacks, and thus we have
$|{\mathbb M}^I|=2^2$, $|{\mathbb M}^A|=2^6$, and $|{\mathbb M}|=2^{6} \times 2^2=256$.
%we have ${\mathbb M} = \{\dbinom{\mathbb N}{1},\cdots,\dbinom{\mathbb N}{10}\}$ where ${\mathbb N} = \{28, 30, 32, 112, 114, 118, 120, 196, 259, 283\}$ is a set of all the possible attack locations.
We, however, could reduce the number of modes by lines 3-6 in Algorithm~\ref{algo5} with ${\TT}=\infty$ into four; i.e., ${\mathbb M}^D_{[1,\infty]}=\{j_1,j_2,j_3,j_4\}$, where each mode associates with one $j \in {\mathbb M}^I$.
%To obtain the reduced mode set, we linearized the nonlinear coupling into $P_{il}^{1}(t) \simeq t_{il}(\theta_i(t)-\theta_l(t))$.
All the four modes assume that sensors (electrical power outputs) $26, 53, 54$ and actuators $14,15$, and $16$ are attacked. Their assumptions on line switching attacks are described as below:\\
Mode $j_1$: There is no line switching attack.\\
Mode $j_2$: Line switch $\{59,60\}$ is attacked.\\
Mode $j_3$: Line switch $\{56,57\}$ is attacked.\\
Mode $j_4$: Line switches $\{56,57\}$, and $\{59,60\}$ are attacked.

Note that the systems under the new modes satisfy Assumptions~\ref{CD_asm_m},~\ref{asm_1},~\ref{asm_2} and uniform observability condition as well as the rest of the assumptions in Theorem~\ref{CD_stability_the1}.
True mode estimation is essential because none of the above modes is true. 
%To find the reduced mode set, we approximate the nonlinear coupling into $P_{il}^{1}(t) \simeq t_{il}(\theta_i(t)-\theta_l(t))$.

We conduct the simulation for Algorithm~\ref{algo5}, using confidence levels $\alpha_1 = \alpha_2 = 0.8$ with corresponding $z$-values $z(\alpha_1)=z(\alpha_2)=1.28$.
%The NISME with Algorithm~\ref{algo4} generates the true mode estimate among $1024$ modes at each time. Since we are only interested in modes 0,1,2, and 3 defined in case study 1, we project all the other modes into mode 4 for the presentation purpose; i.e., if mode 4 is chosen, mode estimation is incorrect.
As case 1, the true modes are among modes 0 to 3, defined in case 1. This is unknown to the defender and the defender remains to consider $256$ possible modes in Algorithm~\ref{algo4}. For the presentation purpose, we project the modes other than modes 0 to 3 into mode 4; i.e., if mode 4 is chosen, mode estimation is incorrect.

\begin{figure}[h]
  \centering
  \includegraphics[width = \linewidth]{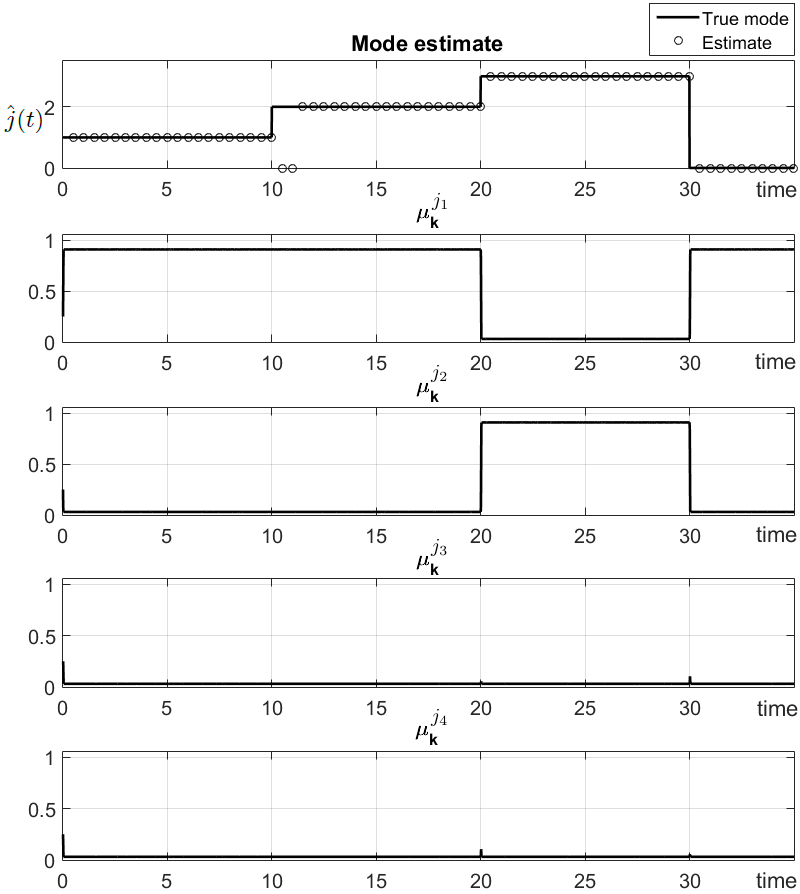}
  \caption{ Mode estimates and probabilities of each mode with reduced mode set.}\label{Case2_2}
\end{figure}

\begin{figure}[h]
  \centering
  \includegraphics[width = \linewidth]{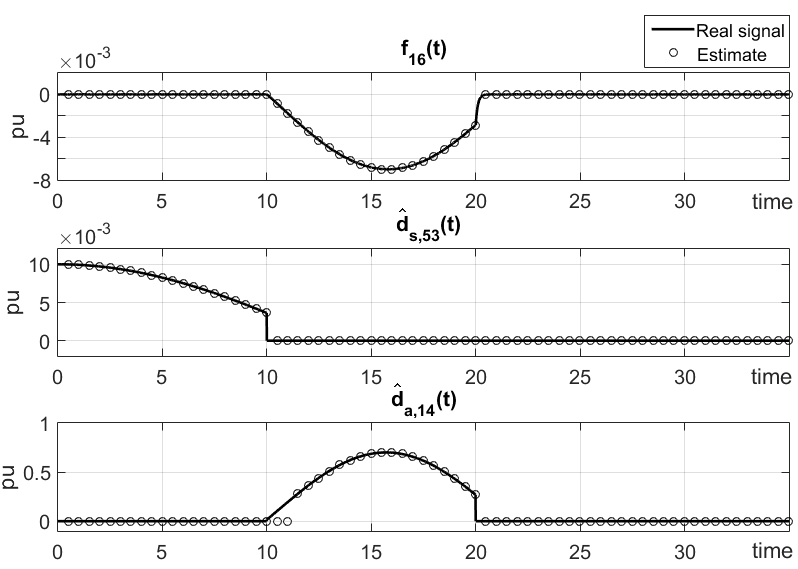}
  \caption{Performance of the NUISE for the reduced mode set. Detailed descriptions on the other subfigures are identical to those of Figure~\ref{Case1_1}. }\label{Case2_1}
\end{figure}

The estimation results are shown in Figure~\ref{Case2_2}, and~\ref{Case2_1}, which are consistent with the results of the case 1 shown in Figure~\ref{Case1_2}, and~\ref{Case1_1}. 
The first subfigure in Figure~\ref{Case2_2} provides a true mode estimation described in Section~\ref{sec:DiscR}.
As case 1, mode estimates are erroneous near $10$ sec because the sizes of attack vectors are small and thus the attack vectors are not regarded statistically significant.
After $30$ sec, mode probabilities oscillate between two modes in case study 1, but not in case study 2.
In case study 1, two modes 1 and 2 are true with zero signal attacks, but mode 3 cannot be a true mode.
In case study 2, only mode $j_1$ is true with zero signal attacks, but modes $j_2$, $j_3$ or $j_4$ cannot be a true mode.

% Sometimes, parts of the attacks are determined statistically insignificant and this results in an incorrect mode estimation.

The modes in the reduced mode set have less restrictive assumptions on attack locations than those of the original mode set, but shows similar estimation results. This simulation, thus, validates the performance of the NISME for the minimal number of modes discussed in Section~\ref{sec:DiscR}.

\section{Estimator derivation and proofs}\label{sec:proofs}
Section~\ref{CD_filter} and~\ref{sec:MM} derive the NISE and and mode estimator, respectively.
The proofs of Theorem~\ref{CD_stability_the1} and Lemma~\ref{Bounded_2} are detailed in Section~\ref{AP_CD_an}.
%Section~\ref{ap_disc} presents the proof of Lemma~\ref{lemm002}.

\subsection{Derivation of the NISE}\label{CD_filter}
We derive the NISE in this section. Since each NISE works only for the corresponding fixed mode, we omit the mode index $j(t)$ for notational simplicity.
To estimate attack vectors and states, we leverage the discrepancies between outputs $z_{1,k}$, $z_{2,k}$, and $z_{3,k}$ defined in~\eqref{E003}, and predicted outputs.
We consider linearization error terms~\eqref{Non_er} in order to establish exact equality, but they would not be considered for the purpose of estimations and covariance updates since they are unknown.
As mentioned before, we derive the NISE by approximating that the attack vector estimates are constants during a sampling period. We denote $\tilde{x}_{k|k,p}$, $\tilde{x}_{k|k-1,p}$, $\tilde{d}_{1,k,p}$, and $\tilde{d}_{2,k,p}$ as estimation errors when the approximation is correct. Additional errors caused by the approximation will be considered in the stability analysis in Section~\ref{AP_CD_an}.

The following lemma guarantees that the derivation steps presented in this section are well defined.
% Note that a positive definite matrix is nonsingular.
\begin{lemma}
If $P_{0}^x$ is positive semidefinite, then, for all $k\geq1$, $P_{k|k-1}^{x}$ is positive semidefinite, and the following matrices induced by the NISE are positive definite:
$\bar{P}_{k|k-1} = C_{3,k}P_{k|k-1}^xC_{3,k}^T+R_{3,k}$, $\tilde{R}_{1,k}\triangleq C_{1,k}P_{k}^x C_{1,k}^T + R_{1,k}$ (if $Rank(\Sigma) \neq 0$), and 
$\tilde{R}_{2,k} \triangleq C_{2,k}(I+\epsilon A_{k-1}-\epsilon G_{1,k-1}M_{1}C_{1,k-1})P_{k-1}^x(I+\epsilon A_{k-1}-\epsilon G_{1,k-1}M_{1}C_{1,k-1})^TC_{2,k}^T+\epsilon^2C_{2,k}G_{1,k-1}$ $M_{1}R_{1,k-1} M_{1}^TG_{1,k-1}^TC_{2,k}^T+\epsilon^2 C_{2,k}Q_{k-1}C_{2,k}^T+R_{2,k}$ (if $Rank(\bar{\Sigma}_k) \neq 0$).
\label{N_singular}
\end{lemma}
\textbf{PROOF.}
We first prove that $P_{k-1}^x \geq 0$, and $P_{k|k-1}^x \geq 0$ for $k \geq 1$ by induction.
Since $P_{0}^x \geq 0$, $P_{1|0}^x \geq 0$ in line 5 of the NISE.
Assume $P_{k-1}^x \geq 0$, and $P_{k|k-1}^x \geq 0$, then $L_k$ is well defined since $R_{3,k}$ are positive definite.
This implies that $P_{k}^x\geq 0$ in line 8 of the NISE, and $P_{k+1}^x\geq 0$ in line 5 of the NISE.
By induction, we conclude that $P_{k-1}^x \geq 0$, and $P_{k|k-1}^x \geq 0$ for $k \geq 1$.

Pick any $k\geq1$, reminding that $P_{k}^x \geq 0$ and $P_{k|k-1}^x \geq 0$.
It holds that $\bar{P}_{k|k-1}$, $\tilde{R}_{1,k}$, and $\tilde{R}_{2,k}$ are positive definite because $R_{1,k}$ (if $Rank(\Sigma) \neq 0$), $R_{2,k}$ (if $Rank(\bar{\Sigma}_k) \neq 0$), and $R_{3,k}$ are positive definite.
Since the above statement holds for any $k\geq1$, we complete the proof.
%The above statement holds for any $k \geq 0$. Thus, now we will show that $P_{k-1}^x$ for $\forall k \geq 0$ are positive semidefinite to complete the proof. From lines 5 and 8 of the NISE, updated covariance matrix $P_{k}^x \geq 0$ is positive semidefinite if $P_{k-1}^x\geq 0$. By induction $P_{k}^x \geq 0$ for $\forall k \geq 0$ since $P_{0}^x \geq 0$ and thus the statement holds.
\oprocend

\subsubsection{Attack vector $d_{1,k-1}$ estimation}
Since attack vector $d_{1,k-1}$ can be directly measured by $z_{1,k-1}$, it can be estimated by, given $\hat{x}_{k-1|k-1}$,
\begin{align}
&\hat{d}_{1,k-1} = M_{1}(z_{1,k-1}-h_1(\hat{x}_{k-1|k-1},u_{k-1},0,t_{k-1}))\nnum\\
&=M_{1}(C_{1,k-1}\tilde{x}_{k-1|k-1,p}+H_{1}d_{1,k-1}+v_{1,k-1}+\psi_{1,k-1})
\label{CD_d1}
\end{align}
where we linearize $h_1$, considering linearization error $\psi_{1,k-1}$.
Estimate $\hat{d}_{1,k-1}$ is unbiased if $M_{1}H_1=I$.
By normalizing the covariance matrix of the right hand side of~\eqref{CD_d1}, we can apply Gauss Markov theorem to find the optimal gain, where best linear unbiased estimate (BLUE) is the linear estimate with smallest variance among all linear and unbiased estimates.
\begin{theorem} (Gauss Markov Theorem~\cite{kailath2000linear})
Estimate $\hat{x} = (H^*H)^{-1}H^*y$ is the BLUE for the model $y = H x + v$ where $v$ is a zero-mean random variable with unit variance and $H$ has full column rank. 
\label{theorem_GM}
\end{theorem}
The BLUE gain can be obtained through Gauss Markov theorem:
\begin{align*}
M_{1}  = (H_{1}^T \tilde{R}_{1,k-1}^{-1}H_{1})^{-1}H_{1}^T\tilde{R}_{1,k-1}^{-1}=H_{1}^{-1}
\end{align*}
because $\tilde{R}_{1,k-1}$ is nonsingular by Lemma~\ref{N_singular} if $Rank(\Sigma) \neq 0$.
Error dynamics for $\tilde{d}_{1,k-1,p}$ is
\begin{align}
&\tilde{d}_{1,k-1,p} =-M_{1}(C_{1,k-1}\tilde{x}_{k-1|k-1,p}+v_{1,k-1}+\psi_{1,k-1}).
\label{error_d1}
\end{align}

\subsubsection{Attack vector $d_{2,k-1}$ estimation}\label{CD_sec_d2}
Attack vector $d_{2,k-1}$ can be estimated using the discrepancy between the predicted output and the measured output
\begin{align}
&\hat{d}_{2,k-1} = M_{2,k} (z_{2,k}-h_2(u_k,0,t_k)\nnum\\
&-C_{2,k}(\hat{x}_{k-1|k-1}+\epsilon f(\hat{x}_{k-1|k-1},u_{k-1},\hat{d}_{1,k-1},0,t_{k-1})))\nnum\\
&=M_{2,k} (C_{2,k}(I+\epsilon A_{k-1})\tilde{x}_{k-1|k-1,p}+\epsilon C_{2,k}\rho_{k-1}\nnum\\
&+\epsilon C_{2,k}G_{1,k-1}\tilde{d}_{1,k-1,p}+\epsilon C_{2,k}G_{2,k-1}d_{2,k-1}\nnum\\
&+\epsilon C_{2,k}w_{k-1}+\epsilon C_{2,k}\phi_{k-1}+v_{2,k}+\psi_{2,k})
\label{CD_d2}
\end{align}
where we linearize functions $h_2$ and $f$ with considering linearization errors~\eqref{Non_er}.
Estimate $\hat{d}_{2,k-1}$ is unbiased if $\epsilon M_{2,k}C_{2,k}G_{2,k-1}=I$.
By Gauss Markov theorem, BLUE gain $M_{2,k}$ can be found by
\begin{align*}
&M_{2,k}=(\epsilon^2 G_{2,k-1}^TC_{2,k}^T\tilde{R}_{2,k}^{-1}C_{2,k}G_{2,k-1})^{-1}\nnum\\
&\times \epsilon G_{2,k-1}^TC_{2,k}^T\tilde{R}_{2,k}^{-1}=(\epsilon C_{2,k}G_{2,k-1})^{-1}
\end{align*}
because $C_{2,k}G_{2,k-1}=C_{2,k}G_{k-1}V_{2}=\bar{\Sigma}_k$ is nonsingular and $\tilde{R}_{2,k}^{-1}$ is also nonsingular as shown in Lemma~\ref{N_singular}. Estimation error is
\begin{align}
&\tilde{d}_{2,k-1,p}=-M_{2,k} (C_{2,k}(I+\epsilon A_{k-1})\tilde{x}_{k-1|k-1,p}\nnum\\
&+\epsilon C_{2,k}G_{1,k-1}\tilde{d}_{1,k-1,p}+\epsilon C_{2,k}w_{k-1}+\epsilon C_{2,k}\phi_{k-1}+v_{2,k}\nnum\\
&+\psi_{2,k}+\epsilon C_{2,k}\rho_{k-1}).
\label{error_d2}
\end{align}

\subsubsection{State prediction}
%Instead, we use the following discretized model in the prediction step \begin{align*} \hat{x}_{k|k-1} &= \hat{x}_{k-1|k-1}+\epsilon f(\hat{x}_{k-1|k-1},u_{k-1},0,t_{k-1})+\epsilon G_{k-1}\hat{d}_{k-1} \end{align*}
Since we have the state and attack vector estimates of the previous time instant, we can predict the current state using the dynamic
\begin{align}
\dot{\hat{x}}(t) = f(\hat{x}(t),u(t),\hat{d}_{1,k-1},0,t)+G_{2,k-1}\hat{d}_{2,k-1}
\label{CD005}
\end{align}
for $t \in (t_{k-1},t_k]$ with initial condition $\hat{x}_{k-1|k-1}$ to have $\hat{x}_{k|k-1}$ at $t=t_k$.
Or equivalently,
\begin{align*}
\hat{x}_{k|k-1} &= \hat{x}_{k-1|k-1}+\epsilon f(\hat{x}_{k-1|k-1},u_{k-1},\hat{d}_{1,k-1},0,t_{k-1})\nnum\\
&+ \epsilon G_{2,k-1}\hat{d}_{2,k-1}+\epsilon\rho_{k-1}
%&=\hat{x}_{k-1|k-1}+\epsilon (A_{k-1} \hat{x}_{k-1|k-1} + B_{k-1}u_{k-1}+ G_{k-1}\hat{d}_{k-1})
\end{align*}
where discretization error $\epsilon\rho_{k-1}$ is considered to obtain the equivalence of~\eqref{CD005}.
Its error dynamic is obtained by applying dynamics~\eqref{DD001} to $x_k$:
\begin{align}
&\tilde{x}_{k|k-1,p} = x_k - \hat{x}_{k|k-1}= (I+\epsilon A_{k-1})\tilde{x}_{k-1|k-1,p}\nnum\\ &+\epsilon G_{1,k-1}\tilde{d}_{1,k-1,p}+\epsilon G_{2,k-1}\tilde{d}_{2,k-1,p}+\epsilon w_{k-1}+\epsilon \phi_{k-1}
\label{p_xkm}
\end{align}
where we linearize $f$ and consider linearization error $\phi_{k-1}$.
The corresponding error covariance matrix is obtained by $P_{k|k-1}^{x}=\bar{A}_{k-1}P_{k-1}^{x}(\bar{A}_{k-1})^T+\bar{Q}_{k-1}$ where we apply~\eqref{error_d1} and~\eqref{error_d2}.

%with the corresponding covariances \begin{align} &P_{k|k-1}^x = (I+\epsilon A_{k-1})P_{k-1}^x(I+\epsilon A_{k-1})^T+Q_{k-1}\nnum\\ &+\epsilon^2 G_{1,k-1}P_{k-1}^{d_1}G_{1,k-1}^T+\epsilon (I+\epsilon A_{k-1})P_{k-1}^{xd_1}G_{1,k-1}^T\nnum\\ &+\epsilon G_{1,k-1}P_{k-1}^{d_1x}(I+\epsilon A_{k-1})^T+\epsilon G_{2,k-1}P_{k-1}^{d_2x}(I+\epsilon A_{k-1})^T\nnum\\ &+\epsilon^2 G_{2,k-1}P_{k-1}^{d_2}G_{2,k-1}^T+\epsilon (I+\epsilon A_{k-1})P_{k-1}^{xd_2}G_{2,k-1}^T\nnum\\ &+\epsilon^2 G_{1,k-1}P_{k-1}^{d_1d_2}G_{2,k-1}^T+\epsilon^2 G_{2,k-1}P_{k-1}^{d_2d_1}G_{1,k-1}^T. \label{CD_Px} \end{align}

\subsubsection{State estimation}
We correct the predicted state estimation using the difference between the expected output and measured output $z_{3,k}$ as
\begin{align}
\hat{x}_{k|k} = \hat{x}_{k|k-1}+L_{k}(z_{3,k} -h_3(\hat{x}_{k|k-1},u_k,0,t_k))
\label{pre_xkk}
\end{align}
and its estimation error is
\begin{align}
\tilde{x}_{k|k,p} =(I-L_{k}C_{3,k})\tilde{x}_{k|k-1,p} -L_k v_{3,k}-L_k\psi_{3,k}
\label{p_xkk}
\end{align}
where we linearize $h_3$ and consider linearization errors $\psi_{3,k}$.
Its error covariance matrix is found by
\begin{align*}
P_{k}^x &=(I-L_{k}C_{3,k})P_{k|k-1}^x(I-L_{k}C_{3,k})^T +L_kR_{3,k}L_k^T.
\end{align*}
Minimizing the error covariance $\trace (P_{k}^x)$ with decision variable $L_{k}$ is an unconstrained optimization problem. To find the BLUE, we take $\trace (P_{k}^x)$ derivative and set it equal to zero $\frac{\partial \trace (P_{k}^x)}{\partial L_k}=2((C_{3,k}^jP_{k|k-1}^{x,j}$ $(C_{3,k}^j)^T+R_{3,k}^j)L_k^T-C_{3,k}P_{k|k-1}^x)=0$
where $L_k=P_{k|k-1}^xC_{3,k}^T (C_{3,k}^jP_{k|k-1}^{x,j}(C_{3,k}^j)^T+R_{3,k}^j)^{-1}$ is the solution of the unconstrained optimization problem and $L_k$ is well defined because $P_{k|k-1}^{x,j} \geq 0 $ (by Lemma~\ref{N_singular}), and $R_{3,k}^j >0$.

%The optimal estimate gain $M_{1,k}=H_{1,k}^{-1}$ can be obtained through Gauss Markov theorem.
%where its error dynamic is \begin{align*} \tilde{d}_{1,k} &=-M_{1,k}(C_{1,k}\tilde{x}_{k|k}+v_{1,k}+\psi_{1,k}). \end{align*} The optimal estimate gain $M_{1,k}=H_{1,k}^{-1}$ can be obtained by following the same argument of Section~\ref{CD_sec_d2}. Corresponding covariance matrices are \begin{align} P_{k}^{d_1} &=  M_{1,k}C_{1,k}P_{k}^xC_{1,k}^TM_{1,k}^T+M_{1,k}R_{1,k}M_{1,k}^T\nnum\\ P_k^{d_1x} &=-M_{1,k}C_{1,k}P_{k}^x. \label{CD_Pd1} \end{align}

Plugging~\eqref{error_d1},~\eqref{error_d2}, and~\eqref{p_xkm} into~\eqref{p_xkk} yields the following update rule for the state estimation errors of the NISE:
\begin{align}
&\tilde{x}_{k|k,p} %(I-L_k C_{3,k}) ((I+\epsilon A_{k-1})\tilde{x}_{k-1|k-1}+\epsilon G_{1,k-1}\tilde{d}_{1,k-1}\nnum\\ &+\epsilon G_{2,k-1}\tilde{d}_{2,k-1}+w_{k-1}+\epsilon \phi_{k-1}) - L_k v_{3,k} -L_k \psi_{k}\nnum\\
%&= (I-L_k C_{3,k}) ((I+\epsilon A_{k-1})\tilde{x}_{k-1|k-1}+\epsilon G_{1,k-1}\tilde{d}_{1,k-1}\nnum\\&-\epsilon G_{2,k-1}M_{2,k} (C_{2,k}(I+\epsilon A_{k-1})\tilde{x}_{k-1|k-1}+v_{2,k}+\psi_{2,k}\nnum\\&+C_{2,k}\phi_{k-1}+\epsilon C_{2,k}G_{1,k-1}\tilde{d}_{1,k-1}+\epsilon C_{2,k}\rho_{k-1}+C_{2,k}w_{k-1})\nnum\\&+w_{k-1}+\phi_{k-1}) - L_k v_{3,k} -L_k \psi_{k}\nnum\\&= (I-L_k C_{3,k}) ((I-\epsilon G_{2,k-1}M_{2,k} C_{2,k})(I+\epsilon A_{k-1})\tilde{x}_{k-1|k-1}\nnum\\&+\epsilon (I-\epsilon G_{2,k-1}M_{2,k} C_{2,k} ) G_{1,k-1}\tilde{d}_{1,k-1}\nnum\\&-\epsilon G_{2,k-1}M_{2,k} v_{2,k}-\epsilon G_{2,k-1}M_{2,k}\psi_{2,k}-\epsilon^2 G_{2,k-1}M_{2,k}C_{2,k}\rho_{k-1}\nnum\\&+(I-\epsilon G_{2,k-1}M_{2,k}C_{2,k})w_{k-1}+(I-\epsilon G_{2,k-1}M_{2,k}C_{2,k}\phi_{k-1})\phi_{k-1}) \nnum\\&- L_k v_{3,k} -L_k \psi_{k}\nnum\\&= (I-L_k C_{3,k}) ((I-\epsilon G_{2,k-1}M_{2,k} C_{2,k})(I+\epsilon A_{k-1})\tilde{x}_{k-1|k-1}\nnum\\&-\epsilon (I-\epsilon G_{2,k-1}M_{2,k} C_{2,k} ) G_{1,k-1}M_{1,k-1}(C_{1,k-1}\tilde{x}_{k-1|k-1}+v_{1,k-1}+\psi_{1,k-1})\nnum\\&-\epsilon G_{2,k-1}M_{2,k} v_{2,k}-\epsilon G_{2,k-1}M_{2,k}\psi_{2,k}-\epsilon^2 G_{2,k-1}M_{2,k}C_{2,k}\rho_{k-1}\nnum\\&+(I-\epsilon G_{2,k-1}M_{2,k}C_{2,k})w_{k-1}+(I-\epsilon G_{2,k-1}M_{2,k}C_{2,k}\phi_{k-1})\phi_{k-1}) \nnum\\&- L_k v_{3,k} -L_k \psi_{k}\nnum\\
=\bar{A}_{k-1}\tilde{x}_{k-1|k-1,p}+\bar{w}_{k-1}+\bar{\phi}_{k-1}+\epsilon\bar{\rho}_{k-1}\nnum\\
&-L_k (C_{3,k} (\bar{A}_{k-1}\tilde{x}_{k-1|k-1,p}+\bar{w}_{k-1}+\bar{\phi}_{k-1}+\epsilon\bar{\rho}_{k-1})\nnum\\
&+v_{3,k}+\psi_{3,k})
\label{CD004}
\end{align}
where $\bar{A}_{k-1}$ and $\bar{w}_{k-1}$ are defined in~\eqref{CD004.7},
\begin{align*}
\bar{\phi}&_{k-1} \triangleq \epsilon (I-\epsilon G_{2,k-1}M_{2,k}C_{2,k})\nnum\\
&\times(\phi_{k-1}- G_{1,k-1}M_{1}\psi_{1,k-1})-\epsilon G_{2,k-1}M_{2,k} \psi_{2,k}\nnum\\
\bar{\rho}&_{k-1} \triangleq -\epsilon G_{2,k-1}M_{2,k} C_{2,k}\rho_{k-1}.
\end{align*}

\subsubsection{The priori probability of the mode}
The priori probability of the mode is derived and explained in the following section. 

\subsection{Derivation of the mode estimator}\label{sec:MM}
%The mode estimator bases on statistical properties of the estimated states. However, we cannot guarantee a special distribution for them; e.g., Gaussian distribution, due to the nonlinearity of the system dynamics. To overcome the difficulties, we assume that the state prediction error is a Gaussian random vector.
It is natural that the predicted output must be matched with the measured output if the mode $j$ is the true mode.
For $\forall j \in {\mathbb M}$, we quantify the discrepancy between the predicted output and the measured output as follows
\begin{align*}
\nu_k^{j} = z_{3,k}^{j} -h_3(\hat{x}_{k|k-1}^{j},u_k,0,j,t_k).
\end{align*}
%The more accurate mode has the smaller discrepancy as well as the smaller covariance of it.
We approximate the output error $\nu_k^{j}$ as a multivariate Gaussian random variable. Then, the likelihood function is given by
\begin{align*}
{\mathcal N}_k^{j}
&\triangleq {\mathcal P}(y_k|j={\rm true})=
{\mathcal N}(\nu_k^{j};0,\bar{P}_{k|k-1}^{j})\nonumber\\
&=\frac{\exp(-(\nu_k^{j})^T(\bar{P}_{k|k-1}^{j})^{-1}\nu_k^{j}/2)}{(2 \pi)^{n^{j}/2} |\bar{P}_{k|k-1}^{j}|^{1/2}}
\end{align*}
where $\bar{P}_{k|k-1}^j = C_{3,k}^jP_{k|k-1}^j(C_{3,k}^j)^T+R_{3,k}^{j}$ is the error covariance matrix of $\nu_k^j$ and $n^j = Rank(\bar{P}_{k|k-1}^j)$.
The likelihood function is well-defined since $\bar{P}_{k|k-1}^{j}>0$ as shown in Lemma~\ref{N_singular}.
By the Bayes' theorem, the posterior probabilities are
$\mu_k^{j} \triangleq {\mathcal P}(j={\rm true}|y_k,\cdots,y_0) =\frac{{\mathcal P}(y_k|j={\rm true}){\mathcal P}(j={\rm true}|y_{k-1},\cdots,y_0)}{\sum_{i=1}^{|{\mathbb M}|}{\mathcal P}(y_k|i={\rm true}){\mathcal P}(i={\rm true}|y_{k-1},\cdots,y_0)}=\frac{{\mathcal N}_k^j\mu_{k-1}^j}{\sum_{i=1}^{|\mathbb M|}{\mathcal N}_k^i\mu_{k-1}^i}$.
However, such update might allow that some $\mu_k^j$ converge to zero.
To prevent this, we modify the posterior probability update to
$
\mu_k^{j} = \frac{\bar{\mu}_k^{j}}{\sum_{i=1}^{|{\mathbb M}|}\bar{\mu}_k^i},
$
where $\bar{\mu}_k^j = \max\{\frac{{\NN}_k^j \mu_{k-1}^j}{\sum_{i=1}^{|{\mathbb M}|}{\NN}_k^i \mu_{k-1}^i},\delta\}$ and
$\delta>0$ is a pre-selected small constant preventing the vanishment of the mode probabilities. 
The last step is to generate the state, attack vector, and mode estimates of the mode having the maximum posteriori probability.

\subsection{Stability analysis of the NISE}\label{AP_CD_an}
The piece-wise constant approximation of the attack vector estimates affects the state estimation error update rule; i.e.,~\eqref{CD004} would be the exact expression only if the approximation is correct.
We first obtain the exact relation between $\tilde{x}_{k|k}$ and $\tilde{x}_{k-1|k-1}$, by analyzing the approximation errors.

Without loss of generality, if $Rank(\Sigma) = 0$ in Assumption~\ref{asm_1},
consider constants $\bar{c}_{1}=\bar{g}_1=\bar{m}_1=0$ for this section.
Likewise, if $Rank(\bar{\Sigma}_k) =0$ for $\forall k$, consider constants $\bar{c}_{2}=\underline{g}_2=\bar{g}_2=\underline{m}_2=\bar{m}_2=\underline{r}_2=0$ for this section.

%We first study the effect of the piece-wise constant approximation for the attack vectors to find the exact relation between $\tilde{x}_{k|k}$ and $\tilde{x}_{k-1|k-1}$.

\subsubsection{Approximation error analysis}\label{CD_error_boundss}
Subscript $e$ is adopted to express \textit{the additional error} induced by the piece-wise constant approximation; e.g., $\tilde{d}_1(t)=d_1(t)-\hat{d}_{1,k}=\tilde{d}_{1,k,p}+d_{1,e}(t)$ for $t \in [t_k,t_{k+1})$.

Under Assumption~\ref{CD_asm1}, the approximation errors for the attack vector estimates, state estimates and state predictions are bounded.
\begin{lemma}
Consider the NISE.
Given Assumptions~\ref{CD_asm1},~\ref{asm_1}, and~\ref{asm_11}, the approximation errors for the attack vector estimates, state estimates and state predictions are bounded by
\begin{align*}
&\|d_{1,e}(t)\| \leq \epsilon\bar{d}, 
&&\|d_{2,e}(t)\| \leq \epsilon\bar{d},\nnum\\
&\|x_e(t)\| \leq \bar{\chi}, &&\|x_{k|k,e}\| \leq (1+\bar{l}\bar{c}_3)\bar{\chi}
\end{align*}
where $\bar{\chi} \triangleq \bar{\Phi}(
\frac{1}{2}\bar{g}_{1}+\bar{g}_{2})\bar{d}\epsilon^2$, $\bar{l} \triangleq \bar{p}\bar{c}_3 \underline{r}_3$,\\
%$\bar{g}_{1} \triangleq \sup_{t}\|G_1(t)\|$, $\bar{g}_{2} \triangleq \sup_{t}\|G_2(t)\|$,
$\bar{\Phi} \triangleq \sup_k\sup_{t,\tau \in [t_k,t_{k+1})}\|\Phi(t,\tau)\|$ and state transition matrix $\Phi(t,\tau) = e^{\int_\tau^t A(\sigma) d \sigma}$.
%for $t \in [t_k,t_{k+1})$, and $t \in (t_k,t_{k+1}]$, respectively.
\label{CD_L_B1}
\end{lemma}
\textbf{PROOF.} 
We prove bounds for $d_{1,e}(t)$ and $d_{2,e}(t)$.
Note that $d_{1,e}(t)$ and $d_{2,e}(t)$ satisfy $\|d_{1,e}(t_k)\|$ $=\|d_{2,e}(t_{k})\|=0$ for $\forall k$
after a new estimate is made at the sampling instants.
Therefore, we consider their bounds during the one sampling time interval only.
Its bound can be found by
\begin{align*}
\|d_{1,e}(t_1)\| &\leq (t_1-t_k) \sup_{t \in [t_k,t_{k+1})}\|\frac{d(t)-d(t_k)}{t-t_k}\| \leq \epsilon\bar{d}
\end{align*}
where $t_1-t_k \leq \epsilon$ is applied to get the desired result.
The proof for $d_{2}(t_2)$ is analogous to that of $d_{1}(t_1)$.
The approximation errors for the attack vector estimates induce additional errors for the state predictions.

We proceed to prove bounds for $x_e(t)$.
Consider the continuous time state prediction error
\begin{align*}
\dot{\tilde{x}}(t)&=f(x(t),u(t),d_1(t),w'(t),t)+G_2(t)d_{2}(t)\nnum\\
&-f(\hat{x}(t),u(t),\hat{d}_{1,k},0,t)-G_{2}(t)\hat{d}_{2,k}\nnum\\
&=A(t)(\tilde{x}_p(t)+x_e(t))+G_{1}(t)(\tilde{d}_{1,k}+d_{1,e}(t))\nnum\\
&+G_{2}(t)(\tilde{d}_{2,k}+d_{2,e}(t))+\epsilon w(t)+\phi(t)
\end{align*}
where $\tilde{x}_p(t)$ represents the state prediction error by prediction~\eqref{CD005} when the approximation is correct, and $x_e(t)=0$ if the approximation is correct.
Therefore, we consider the bound during the one sampling interval and it holds that
\begin{align*}
\dot{x}_e(t) &=A(t)x_e(t)+G_{1}(t)d_{1,e}(t)+G_{2}(t)d_{2,e}(t).
\end{align*}
For $t \in [t_k,t_{k+1}]$, its solution is
\begin{align*}
x_e(t) & = \int_{t_k}^{t_{k+1}} \Phi(t,\tau)(G_{1}(\tau)d_{1,e}(\tau)+G_{2}(\tau)d_{2,e}(\tau))d \tau.
\end{align*}
Its bound can be found by, taking norm both sides,
\begin{align*}
&\|x_e(t)\| \leq \int_{t_k}^t \bar{\Phi}(\bar{g}_{1}(\tau-t_{k})\bar{d}+\bar{g}_{2}(t_{k+1}-\tau)\bar{d} \ d \tau \nnum\\
&=\bar{\Phi}(
\frac{1}{2}\bar{g}_{1}\bar{d}(t-t_k)+\frac{1}{2}\bar{g}_{2}\bar{d}(t_{k+1}-t)+\frac{1}{2}\bar{g}_{2}\bar{d}\epsilon)(t-t_{k}).
\end{align*}
We can obtain the desired result by applying $t_{k+1}-t \leq \epsilon$ and $t-t_{k} \leq \epsilon$.
From the result, we can derive the approximation error for the state estimations $x_{k|k,e}$.

Consider the state estimation error
\begin{align}
&\tilde{x}_{k|k,p}+x_{k|k,e}=\tilde{x}_{k|k}=x_k-\hat{x}_{k|k} \nnum\\
&=\tilde{x}_{k|k-1,p}+x_{e}(t_k)-L_{k}(C_{3,k}(\tilde{x}_{k|k-1,p}+x_{e}(t_k))\nnum\\
&+v_{3,k}+\psi_{3,k})
\label{CD009}
\end{align}
where~\eqref{pre_xkk} is applied to the third equality.
Subtracting $\tilde{x}_{k|k,p}$ in~\eqref{CD004} both sides, we have
\begin{align}
x_{k|k,e}&=x_{e}(t_k)-L_{k}C_{3,k}x_{e}(t_k).
\label{CD010}
\end{align}
By taking the norm both sides and apply Cauchy-Schwarz inequality, its upper bound can be found by $\|\tilde{x}_{k|k,e}\| \leq (1+\bar{l}\bar{c}_3)\bar{\chi}$.
\oprocend

\subsubsection{Proof of Theorem~\ref{CD_stability_the1}}\label{AP_CD_Th_proof}
We define Lyapunov candidate $V_{k} \triangleq \tilde{x}_{k|k}^T (P_{k}^x)^{-1}\tilde{x}_{k|k}$ with $0<\underline{p}I \leq P_k^x \leq \bar{p}I$.
Plugging~\eqref{CD004} and~\eqref{CD010} into the update rule of $\tilde{x}_{k|k}$ in~\eqref{CD009}, we obtain
\begin{align}
&V_{k}%=((I-L_{k}C_{3,k})(\bar{A}_{k-1}\tilde{x}_{k-1|k-1}+x_{e}(t_{k}))+\bar{v}_k\nnum\\&+\bar{\psi}_k-L_k\psi_{k,e})^TP_{k|k}^{-1} ((I-L_{k}C_{3,k})(\bar{A}_{k-1}\tilde{x}_{k-1|k-1}+x_{e}(t_{k}))\nnum\\&+\bar{v}_k+\bar{\psi}_k-L_k\psi_{k,e})\nnum\\&
=\tilde{x}_{k-1|k-1}^T\bar{A}_{k-1}^T(I-L_{k}C_{3,k})^T(P_{k}^x)^{-1}(I-L_{k}C_{3,k})\nnum\\
&\times \bar{A}_{k-1}\tilde{x}_{k-1|k-1}
+\bar{v}_k^T(P_{k}^x)^{-1}\nnum\\
&\times (2(I-L_{k}C_{3,k})(\bar{A}_{k-1}\tilde{x}_{k-1|k-1}+x_e(t_{k}))+\bar{v}_k+2\bar{\psi}_k)\nnum\\
&+\bar{\psi}_k^T(P_{k}^x)^{-1}(2(I-L_{k}C_{3,k})\bar{A}_{k-1}\tilde{x}_{k-1|k-1}+\bar{\psi}_k)\nnum\\
&+x_{e}^T(t_{k})(I-L_{k}C_{3,k})^T(P_{k}^x)^{-1}\nnum\\
&\times((I-L_{k}C_{3,k})(2\bar{A}_{k-1}\tilde{x}_{k-1|k-1}+x_{e}(t_{k}))+2\bar{\psi}_k)
\label{CD013}
\end{align}
where $\bar{\psi}_k \triangleq (I-L_{k}C_{3,k})(\bar{\phi}_{k-1}+\epsilon\bar{\rho}_{k-1})-L_{k}\psi_{3,k}$
and $\bar{v}_k \triangleq (I-L_{k}C_{3,k})\bar{w}_{k-1}-L_{k}v_{3,k}$.
To deal with the above terms, we formalize the proof by several claim statements; i.e., Claim B.1-B.4 deal with the first to fourth terms in~\eqref{CD013}, respectively.

Recall that $Q_k \leq \bar{q}' I$, $R_{1,k} \leq \bar{r}_1$, $R_{2,k} \leq \bar{r}_2$, $R_{3,k} \leq \bar{r}_3$, and $\epsilon \leq \bar{\epsilon}$ for some constants $\bar{q}'$, $\bar{r}_1$, $\bar{r}_2$, $\bar{r}_3$, and $\bar{\epsilon}$. We will show that, for any $\gamma \in (0,1)$, if $\bar{q}'$, $\bar{r}_1$, $\bar{r}_2$, $\bar{r}_3$, and $\bar{\epsilon}$ are chosen properly, then the state estimation errors $\tilde{x}_{k|k}$ and attack vector estimation errors $\tilde{d}_1(t)$, $\tilde{d}_2(t)$ hold PESp-like properties.
Assumption~\ref{asm_1} implies that
$\|\bar{A}_{k} \|\leq \bar{a} \triangleq (1+\epsilon\bar{g}_2\bar{c}_2\bar{m}_2)(1+\epsilon\bar{a}'+\epsilon\bar{g}_1\bar{c}_1\bar{m}_1)$, and
$\underline{q} I \leq \bar{Q}_k \leq \bar{q} I$ where 
$\underline{q} \triangleq \max\{\epsilon^2\underline{g}_2^2\underline{m}_2^2\underline{r}_2,\epsilon^2\underline{q}'\}$ and
$\bar{q} \triangleq \epsilon^2(1+\epsilon\bar{g}_2\bar{c}_2 \bar{m}_2)^2(\bar{q}'+\bar{g}_1\bar{m}_1\bar{r}_1)^2+\epsilon^2\bar{g}_2^2\bar{m}_2^2\bar{r}_2$.

\textbf{Claim B.1: }
There exists a constant $\alpha' = (1+\frac{\underline{q}}{2\bar{a}^2\bar{p}})^{-1} \in(0,1)$ such that
\begin{align*}
&\bar{A}_{k-1}^T(I-L_{k}C_{3,k})^T(P_{k}^x)^{-1}(I-L_{k}C_{3,k})\bar{A}_{k-1}\nnum\\
&<\alpha' (P_{k-1}^{x})^{-1}.
\end{align*}
\textbf{PROOF.} 
Error covariance matrix is updated as
\begin{align*}
P_{k}^x &= (I-L_{k}C_{3,k})\bar{A}_{k-1}P_{k-1}^x\bar{A}_{k-1}^T(I-L_{k}C_{3,k})^T\nnum\\
&+(I-L_{k}C_{3,k})\bar{Q}_{k-1}(I-L_{k}C_{3,k})^T+L_kR_{3,k}L_k^T.
\end{align*}
%where $L_k=P_{k|k-1}^xC_{3,k}^T \tilde{R}_{3,k}^{-1}$ with $\tilde{R}_{3,k} \triangleq C_{3,k}P_{k|k-1}^xC_{3,k}^T+R_{3,k}$. 
From the upper and lower bounds of matrices, we have 
\begin{align*}
P_{k}^x &> (I-L_{k}C_{3,k})\bar{A}_{k-1}P_{k-1}^x\bar{A}_{k-1}^T(I-L_{k}C_{3,k})^T\nnum\\
&\times(1+\frac{\underline{q}}{2\bar{a}^2\bar{p}})+L_kR_{3,k}L_k^T.
\end{align*}
This implies that $P_{k}^x-(I-L_{k}C_{3,k})\bar{A}_{k-1}P_{k-1}^x\bar{A}_{k-1}^T(I-L_{k}C_{3,k})^T(1+\frac{\underline{q}}{2\bar{a}^2\bar{p}})>0$ and thus it holds that
\begin{align*}
&(1+\frac{\underline{q}}{2\bar{a}^2\bar{p}}) P_{k-1}^x+(1+\frac{\underline{q}}{2\bar{a}^2\bar{p}})P_{k-1}^x \bar{A}_{k-1}^T(I-L_{k}C_{3,k})^T\nnum\\
&\times (P_{k}^x-(I-L_{k}C_{3,k})\bar{A}_{k-1}P_{k-1}^x\bar{A}_{k-1}^T(I-L_{k}C_{3,k})^T)\nnum\\
&>0.
\end{align*}
By applying the matrix inversion lemma~\cite{hager1989updating}, it follows that
\begin{align*}
&((1+\frac{\underline{q}}{2\bar{a}^2\bar{p}})^{-1}(P_{k-1}^x)^{-1}-\bar{A}_{k-1}^T(I-L_{k}C_{3,k})^T(P_{k}^x)^{-1}\nnum\\
&\times (I-L_{k}C_{3,k})\bar{A}_{k-1})^{-1}>0
\end{align*}
which implies the statement with $\alpha' = (1+\frac{\underline{q}}{2\bar{a}^2\bar{p}})^{-1}$.
\oprocend
%By applying Lemma 6.1 in~\cite{kluge2010stochastic}, we have the desired result.

\textbf{Claim B.2:}
There exists a positive constant $\epsilon_0>0$ such that
\begin{align*}
&{\mathbb E}[\bar{v}_k^T(P_{k}^x)^{-1}(2(I-L_{k}C_{3,k})(\bar{A}_{k-1}\tilde{x}_{k-1|k-1}+x_e(t_{k}))\nnum\\
&+\bar{v}_k+2\bar{\psi}_k)] \leq \epsilon_0.
\end{align*}
\textbf{PROOF.}
Noises $w_{k-1}, v_{1,k-1}, v_{2,k}$, and $v_{3,k}$ are uncorrelated and thus we have
\begin{align*}
&{\mathbb E}[\bar{v}_k^T(P_{k}^x)^{-1}(2(I-L_{k}C_{3,k})(\bar{A}_{k-1}\tilde{x}_{k-1|k-1}+x_e(t_{k}))\nnum\\
&+\bar{v}_k+2\bar{\psi}_k)]\nnum\\
&={\mathbb E}[\bar{v}_k^T(P_{k}^x)^{-1}\bar{v}_k] \leq \epsilon^2\underline{p}(1+\bar{l}\bar{c}_3)^2((1+\epsilon \bar{g}_2 \bar{m}_2 \bar{c}_2)^2\nnum\\
&\times (\bar{q}{Rank(Q_{k-1})}+ \bar{g}_1^2\bar{m}_1^2\bar{r}_1 {Rank(R_{1,k-1})})\nnum\\
&+\bar{g}_2^2\bar{m}_2^2\bar{r}_2{Rank(R_{2,k})})+\underline{p}\bar{l}^2\bar{r}_3{Rank(R_{3,k})} \triangleq \epsilon_0
\end{align*}
where we apply $\|v_{1,k}\|^2 =\trace (v_{1,k}^T v_{1,k})=\trace (v_{1,k} v_{1,k}^T) \leq \bar{r}_1 {Rank(R_{1,k})}$ and the similar relations for $\|w_{k-1}\|^2$, $\|v_{2,k}\|^2$, and $\|v_{3,k}\|^2$.
\oprocend

\textbf{Claim B.3:}
There exist constants $\delta, \delta_{\rho}, \lambda_1, \epsilon_1>0$ such that, for $\forall \|\tilde{x}_{k-1|k-1}\| \leq \delta$ and $\epsilon \leq \delta_{\rho}$, the following holds:
\begin{align*}
&\bar{\psi}_k^T(P_{k}^x)^{-1}(2(I-L_{k}C_{3,k})\bar{A}_{k-1}\tilde{x}_{k-1|k-1}+\bar{\psi}_k) \nnum\\
&\leq \lambda_1\|\tilde{x}_{k-1|k-1}\|^3+\epsilon_1.
\end{align*}
\textbf{PROOF.} 
By Assumptions~\ref{asm_1} and~\ref{asm_2}, it holds that
\begin{align*}
&\|\bar{\psi}_k\| =\|(I-L_{k}C_{3,k})(\bar{\phi}_{k-1}+\epsilon\bar{\rho}_{k-1})-L_{k}\psi_{3,k}\|\nnum\\
&\leq (1+\bar{l}\bar{c}_3)(\epsilon(1+\bar{g}_2\bar{c}_2\bar{m}_2)(\epsilon_{\phi}+\epsilon_{\psi_1}\bar{g}_1\bar{m}_1)+\epsilon \bar{g}_2\bar{m}_2\epsilon_{\psi_2}) \nnum\\
&\times \|\tilde{x}_{k-1|k-1}\|^2+ \epsilon^2 \epsilon_{\rho}\bar{g}_2\bar{m}_2\bar{c}_2 + \bar{l}\epsilon_{\psi_3}\|\tilde{x}_{k-1|k-1}\|^2\nnum\\
&\triangleq \lambda_1' \|\tilde{x}_{k-1|k-1}\|^2+\epsilon_1'
\end{align*}
for all $\|\tilde{x}_{k-1|k-1}\| \leq \delta$, and $\epsilon \leq \delta_{\rho}$.
Therefore, we have
\begin{align*}
&\bar{\psi}_k^T(P_{k}^x)^{-1}(2(I-L_{k}C_{3,k})\bar{A}_{k-1}\tilde{x}_{k-1|k-1}+\bar{\psi}_k)\nnum\\
&\leq 2\lambda_1' \underline{p}(1+\bar{l}\bar{c}_3)\bar{a}\|\tilde{x}_{k-1|k-1}\|^3+\lambda_1'^2\underline{p}\delta\|\tilde{x}_{k-1|k-1}\|^3\nnum\\
&+\epsilon_1'\underline{p}(2(1+\bar{l}\bar{c}_3)\bar{a}\delta+\lambda_1'\delta^2+\epsilon_1')\nnum\\
%&\leq k_1' \underline{p}(2(1+\frac{1}{k_1'})(1+\bar{l}\bar{c}_3)\bar{a}+(1+k_1')\delta_1+\epsilon_1')\nnum\\&\times \|\tilde{x}_{k-1|k-1}\|^3+\epsilon_1' \underline{p}(\epsilon_1'+2(1+\bar{l}\bar{c}_3)\bar{a}+k_1')\nnum\\
&\triangleq \lambda_1 \|\tilde{x}_{k-1|k-1}\|^3+\epsilon_1
\end{align*}
where $\|\tilde{x}_{k-1|k-1}\| \leq \delta$ is applied.
\oprocend

\textbf{Claim B.4:}
There exist constants $\delta,\lambda_2,\epsilon_2>0$ such that, for $\forall \|\tilde{x}_{k-1|k-1}\| \leq \delta$, the following holds:
\begin{align*}
&x_{e}^T(t_{k})(I-L_{k}C_{3,k})^T(P_{k}^x)^{-1}((I-L_{k}C_{3,k})(2\bar{A}_{k-1}\nnum\\
&\times \tilde{x}_{k-1|k-1}+x_{e}(t_{k}))+2\bar{\psi}_k)\leq \epsilon_2 + \lambda_2\|\tilde{x}_{k-1|k-1}\|^3.
\end{align*}
\textbf{PROOF.} 
By noting that $\|x_{e}(t_{k})\| \leq \bar{\chi}$, we have
\begin{align*}
&x_{e}^T(t_k)(I-L_{k}C_{3,k})^T(P_{k}^x)^{-1}((I-L_{k}C_{3,k})\nnum\\
&\times (2\bar{A}_{k-1}\tilde{x}_{k-1|k-1}+x_{e}(t_{k}))+2\bar{\psi}_k)\nnum\\
&\leq \underline{p}(1+\bar{l}\bar{c}_3)((1+\bar{l}\bar{c}_3)(2\bar{a}\|\tilde{x}_{k-1|k-1}\|+\bar{\chi})\nnum\\
&+2\lambda_1'\|\tilde{x}_{k-1|k-1}\|^2+2\epsilon_1')\bar{\chi} \leq \epsilon_2 + \lambda_2 \|\tilde{x}_{k-1|k-1}\|^3
\end{align*}
for any $\|\tilde{x}_{k-1|k-1}\| \leq \delta$ where %$\|\tilde{x}_{k-1|k-1}\| \leq 1+\|\tilde{x}_{k-1|k-1}\|^3$ and $\|\tilde{x}_{k-1|k-1}\|^2 \leq 1+\|\tilde{x}_{k-1|k-1}\|^3$ are applied, and 
$\epsilon_2 \triangleq 
\underline{p}(1+\bar{l}\bar{c}_3)((1+\bar{l}\bar{c}_3)(2\bar{a}+\bar{\chi})+2\lambda_1'+2\epsilon_1')\bar{\chi}$
and
$\lambda_2 \triangleq \underline{p}(1+\bar{l}\bar{c}_3)(2(1+\bar{l}\bar{c}_3)\bar{a}+2\lambda_1')\bar{\chi}$.
\oprocend

By applying Claims B.1-B.4 to~\eqref{CD013}, we have
\begin{align*}
{\mathbb E}[V_{k}] &\leq \alpha' {\mathbb E}[\tilde{x}_{k-1|k-1}^TP_{k-1}^{-1}\tilde{x}_{k-1|k-1}]
\nnum\\
&+\lambda{\mathbb E}[\|\tilde{x}_{k-1|k-1}\|^3]+(\epsilon_0+\epsilon_1+\epsilon_2)
\end{align*}
for $\forall \|\tilde{x}_{k-1|k-1}\| \leq \delta$ and $\epsilon \leq \delta_{\rho}$ where $\lambda \triangleq \lambda_1+\lambda_2$.
Note that $\lambda_1$ tends to zero as $\epsilon_{\phi},\epsilon_{\psi_1},\epsilon_{\psi_2},\epsilon_{\psi_3},\epsilon_{\rho}$ tend to zero; $\lambda_2$ tends to zero as $\epsilon = t_k-t_{k-1}$ tends to zero.
By choosing a sufficiently small tuple $(\epsilon_{\phi},\epsilon_{\psi_1},\epsilon_{\psi_3},\epsilon_{\rho},\bar{\epsilon}')$ such that $\lambda\delta<\alpha'\bar{p}^{-1}$, we have
\begin{align}
{\mathbb E}[V_{k}] &\leq \alpha{\mathbb E}[V_{k-1}]+c
\label{E005}
\end{align}
for $\|\tilde{x}_{k-1|k-1}\| \leq \delta$ where $0 < \alpha<1$ and $c \triangleq \epsilon_0+\epsilon_1+\epsilon_2$.
Inequality~\eqref{E005} holds for all $\epsilon \leq \bar{\epsilon}'$.
Remind that $\epsilon_0$ tends to zero as $\bar{q}'$, $\bar{r}_1$, $\bar{r}_2$ and $\bar{r}_3$ tend to zero,
and constants $\epsilon_1$ and $\epsilon_2$ tend to zero as $\epsilon$ tends to zero.
Thus, for any given constant $c'>0$, we can choose sufficiently small tuple $(\bar{q}'$, $\bar{r}_1$, $\bar{r}_2$, $\bar{r}_3$, $\bar{\epsilon})$ such that $c<c'$ holds. 
The following claim shows a PESp-like property for state estimation errors.

\textbf{Claim B.5:}
For any $\gamma \in (0,1)$, there exist positive constants $\alpha_x$, $b_x$, $c_x$, $\underline{\delta}$ and tuple $(\bar{q}'$, $\bar{r}_1$, $\bar{r}_2$, $\bar{r}_3$, $\bar{\epsilon})$
such that, for all $\|\tilde{x}_{0|0}\| \leq \underline{\delta}$, the following holds for all $k \geq 0$:
\begin{align*}
P(\|\tilde{x}_{k|k}\| < \alpha_x e^{-b_x k}\|\tilde{x}_{0|0}\|+c_x) \geq 1-\gamma.
\end{align*}
\textbf{PROOF.} 
Consider any $\gamma \in (0,1)$ and $\gamma_1 <\gamma$.
Then, there exists sufficiently small constant $\underline{\delta}<\delta$  and tuple $(\bar{q}'$, $\bar{r}_1$, $\bar{r}_2$, $\bar{r}_3$, $\bar{\epsilon})$ such that
$\bar{p}\underline{\delta}^2+\frac{c}{1-\alpha}\leq \gamma_1 \underline{p}\delta^2$ holds. Since $V \leq\bar{p}\|\tilde{x}\|^2$, we have 
\begin{align}
V+\frac{c}{1-\alpha}\leq\bar{p}\|\tilde{x}\|^2+\frac{c}{1-\alpha}\leq \gamma_1 \underline{p}\delta^2, \ \forall \|\tilde{x}\| \leq \underline{\delta}.
\label{np_001}
\end{align}
We choose any $\|\tilde{x}_{0|0}\| \leq \underline{\delta}$.
Define the first exit time $\mu \triangleq \inf\{t_k> 0| \|\tilde{x}_{k|k}\|>\delta\}$, and $\mu \wedge k \triangleq \min \{\mu,k\}$ for any $k >0$.
We have
\begin{align}
\underline{p}\delta^2 P(\mu \leq k)&= {\mathbb E}[\underline{p}\delta^2I_{[\mu \leq k]}]
\leq {\mathbb E}[ V_{\mu}I_{[\mu \leq k]}]\nnum\\
&\leq{\mathbb E}[V_{\mu \wedge k}]\leq \alpha^{\mu \wedge k} V_0 + c\sum_{i=0}^{\mu \wedge k-1}\alpha^i
\label{np_002}
\end{align}
where indicator function $I_{[\mu \leq k]}$ satisfies $I_{[\mu \leq k]}=1$ if $\mu \leq k$, otherwise $I_{[\mu \leq k]}=0$.
The first inequality holds because $\underline{p}\delta^2<\underline{p}\|\tilde{x}_{\mu|\mu}\|^2\leq V_{\mu}$.
The third inequality is derived by recursively applying ${\mathbb E}[V_k] \leq \alpha {\mathbb E}[V_{k-1}]+c$.
By~\eqref{np_001}, it follows from~\eqref{np_002} that
\begin{align}
P(\mu \leq k) &\leq \gamma_1\frac{\alpha^{\mu \wedge k} V_0 + c\sum_{i=0}^{\mu \wedge k-1}\alpha^i}{V_0 +\frac{c}{1-\alpha}}\leq \gamma_1.
\label{np_003}
\end{align}
Letting $k \rightarrow \infty$, we also have 
\begin{align*}
P(\|\tilde{x}_{k|k}\|\leq\delta) &\geq 1-\gamma_1.
\end{align*}
Again consider any $k$ with $\mu \wedge k$, and $\gamma_2 = \gamma-\gamma_1$.
%For $\|x\| \geq \frac{c \lambda}{2 \underline{p}}$, we have  \begin{align*} {\mathbb E}[V_{t}] \leq \frac{\alpha}{2}{\mathbb E}[V_{t-1}] \leq (\frac{\alpha}{2})^t V_0. \end{align*} Markov's inequality derives \begin{align*} P(V_t \geq \frac{(\alpha/2)^t V_0}{\gamma_2}) &\leq \gamma_2\frac{ {\mathbb E}[V_t]}{(\alpha/2)^t V_0}\leq \gamma_2. \end{align*} This implies that \begin{align} &P(\|\tilde{x}_{t|t}\| \leq  \sqrt{\frac{\bar{p}}{\gamma_2\underline{p}}} (\frac{\alpha}{2})^{t/2}\|\tilde{x}_0\|)=P(\|\tilde{x}_{t|t}\| \leq \beta(\| tilde{x}_0\|,t))\nnum\\ &\leq 1 - \gamma_2. \label{np_004} \end{align} Now consider $\|x\| \leq \frac{c \lambda}{2 \underline{p}}$
%$P(X<Y) \leq \frac{{\mathbb E}[X]}{Y}$ for nonnegative random variables $X$ and $Y$.
Markov's inequality (p.455~\cite{hoffman1994probability}) derives
\begin{align*}
&P(V_{\mu \wedge k} \geq \frac{\alpha^{\mu \wedge k} V_{0}+\frac{c}{1-\alpha}}{\gamma_2}) \leq \gamma_2\frac{ {\mathbb E}[V_{\mu \wedge k}]}{\alpha^{\mu \wedge k} V_{0}+\frac{c}{1-\alpha}}\nnum\\
&\leq \gamma_2\frac{\alpha^{\mu \wedge k} V_{0}+c\sum_{i=0}^{\mu \wedge k-1}\alpha^i}{\alpha^{\mu \wedge k} V_{0}+\frac{c}{1-\alpha}}\leq \gamma_2.
\end{align*}
%where we apply $c\sum_{i=0}^{\mu \wedge k-1}\alpha^i \leq c\sum_{i=0}^{\infty}\alpha^i = \frac{c}{1-\alpha}$.
Equivalently, $P(V_{\mu \wedge k} < \frac{\alpha^{\mu \wedge k} V_{0}+\frac{c}{1-\alpha}}{\gamma_2})\geq 1-\gamma_2$.
This implies that, by Minkowski inequality, 
\begin{align}
&P(\|\tilde{x}_{\mu \wedge k|\mu \wedge k}\| < \sqrt{\frac{\bar{p}}{\gamma_2\underline{p}}}
\alpha^{(\mu \wedge k)/2}\|\tilde{x}_{0|0}\|+\sqrt{\frac{c}{(1-\alpha)(\underline{p}\gamma_2)}})\nnum\\
&=P(\|\tilde{x}_{\mu \wedge k|\mu \wedge k}\| < \beta_x(\|\tilde{x}_{0|0}\|,\mu \wedge k)+c_x)\geq 1 - \gamma_2
\label{np_004}
\end{align}
where $\beta_x(\|\tilde{x}_{0|0}\|,\mu \wedge k)=\sqrt{\frac{\bar{p}}{\gamma_2\underline{p}}}
\alpha^{(\mu \wedge k)/2}\|\tilde{x}_{0|0}\|$ and $c_x =\sqrt{\frac{c}{(1-\alpha)(\underline{p}\gamma_2)}})$.
By~\eqref{np_003} and~\eqref{np_004}, we can obtain %condition $\mu > k$, as follows
\begin{align*}
&P(\|\tilde{x}_{k|k}\| < \beta_x(\|\tilde{x}_{0|0}\|,k)+c_x)\nnum\\
&=P(\|\tilde{x}_{k|k}\| < \beta_x(\|\tilde{x}_{0|0}\|,k)+c_x|\mu>k)P(\mu>k)\nnum\\
&\quad +P(\|\tilde{x}_{k|k}\| < \beta_x(\|\tilde{x}_{0|0}\|,k)+c_x|\mu \leq k)P(\mu \leq k)\nnum\\
&\geq P(\|\tilde{x}_{\mu \wedge k|\mu \wedge k}\| < \beta_x(\|\tilde{x}_{0|0}\|,\mu \wedge k)+c_x| \mu > k)\nnum\\
&\quad \times P(\mu > k)
\end{align*}
\begin{align}
&= P(\|\tilde{x}_{\mu \wedge k|\mu \wedge k}\| < \beta_x(\|\tilde{x}_{0|0}\|,\mu \wedge k)+c_x)\nnum\\
&\quad - P(\|\tilde{x}_{\mu \wedge k|\mu \wedge k}\| < \beta_x(\|\tilde{x}_{0|0}\|,\mu \wedge k)+c_x| \mu \leq k)\nnum\\
&\quad \times P(\mu \leq k)\nnum\\
&\geq 1 - \gamma_2-\gamma_1 =1- \gamma.
\label{np_006}
\end{align}
%If $\alpha \leq \frac{1}{e}$, then $a_x=\sqrt{\frac{\bar{p}}{\gamma_2\underline{p}}}$, and $b_x=0.5$. Otherwise, $a_x=\sqrt{\frac{\bar{p}}{\gamma_2\underline{p}}}$, and $b_x=-0.5 \ln \alpha>0$.
\oprocend

We now proceed to prove PESp-like property of $\tilde{d}_1(t)$.
Without loss of generality, consider $\gamma$, $\gamma_1$, and $\tilde{x}_{0|0}$ used in Claim B.5. 
Consider mean square errors of $\hat{d}_1(t)$ for $t \in [t_k,t_{k+1})$, with equation~\eqref{error_d1} and $\|d_{1,e}(t)\| \leq \epsilon \bar{d}$ by Lemma~\ref{CD_L_B1}:
\begin{align*}
&{\mathbb E}[\|\tilde{d}_1(\mu \wedge t)\|^2]={\mathbb E}[\|\tilde{d}_{1,\mu \wedge k,p}+d_{1,e}(\mu \wedge t)\|^2]\nnum\\
%={\mathbb E}[(M_{1}(C_{1,k}\tilde{x}_{k|k}+v_{1,k}+\psi_{1,k}))^T\nnum\\ &\quad \times M_{1}(C_{1,k}\tilde{x}_{k|k}+v_{1,k}+\psi_{1,k})]\nnum\\
&\leq \bar{m}_1^2({\mathbb E}[\bar{c}_1^2\|\tilde{x}_{\mu \wedge k|\mu \wedge k,p}\|^2+\bar{c}_1\epsilon \bar{d}\|\tilde{x}_{\mu \wedge k|\mu \wedge k,p}\|\nnum\\
&+\|v_{1,\mu \wedge k}\|^2+\|\psi_{1,\mu \wedge k}\|^2+2\bar{c}_1\|\tilde{x}_{\mu \wedge k|\mu \wedge k,p}\|\|\psi_{1,\mu \wedge k}\|\nnum\\
&+\epsilon \bar{d}\|\psi_{1,\mu \wedge k}\|])+\epsilon^2 \bar{d}^2
\end{align*}
where we apply Cauchy-Schwarz inequality.
Since $\|\tilde{x}_{\mu \wedge k|\mu \wedge k,p}\| \leq \|\tilde{x}_{\mu \wedge k|\mu \wedge k}\|+\|{x}_{\mu \wedge k|\mu \wedge k,e}\| \leq \|\tilde{x}_{\mu \wedge k|\mu \wedge k}\|$ $+ (1+ \bar{l}\bar{c}_2) \bar{\chi}$ and $\|\tilde{x}_{\mu \wedge k|\mu \wedge k}\|\leq \delta$, it follows that
\begin{align*}
&{\mathbb E}[\|\tilde{d}_1(\mu \wedge t)\|^2]\leq \bar{m}_1^2(\bar{c}_1^2+2\epsilon_{\psi_1}\bar{c}_1(\delta+(1+ \bar{l}\bar{c}_2) \bar{\chi})\nnum\\
&+\epsilon_{\psi_1}^2(\delta^2+\epsilon \bar{d})){\mathbb E}[\tilde{x}_{\mu \wedge k|\mu \wedge k}]+\bar{m}_1^2(\bar{c}_1^2(\delta(1+ \bar{l}\bar{c}_2) \bar{\chi}\nnum\\
&+(1+ \bar{l}\bar{c}_2)^2 \bar{\chi}^2)+\bar{c}_1\epsilon \bar{d}(\delta+(1+ \bar{l}\bar{c}_2) \bar{\chi})\nnum\\
&+\bar{r}_1 Rank(R_{1,\mu \wedge k}))+\epsilon^2\bar{d}^2
\end{align*}
where $\|v_{1,\mu \wedge k}\|^2 \leq \bar{r}_1 Rank(R_{1,\mu \wedge k})$, and $\|\psi_{1,\mu \wedge k}\| \leq \epsilon_{\psi_1}\|\tilde{x}_{\mu \wedge k|\mu \wedge k}\|^2$ are applied.
Applying $\epsilon \leq c$ and ${\mathbb E}[\|\tilde{x}_{\mu \wedge k|\mu \wedge k}\|^2]$ $\leq \frac{\alpha^{\mu \wedge k}}{\underline{p}} V_0 + \frac{c}{\underline{p}}\sum_{i=0}^{\mu \wedge k-1}\alpha^i$ obtained in~\eqref{np_002}, 
it follows that
\begin{align*}
&{\mathbb E}[\|\tilde{d}_1(\mu \wedge t)\|^2] \leq \beta_1(\|\tilde{x}_{0|0}\|^2,\mu \wedge k)+ c_1
\end{align*}
where
\begin{align*}
&\beta_1(\|\tilde{x}_{0|0}\|^2,\mu \wedge k)=\bar{m}_1^2(\bar{c}_1^2+2\epsilon_{\psi_1}\bar{c}_1(\delta+(1+ \bar{l}\bar{c}_2) \bar{\chi})\nnum\\
&+\epsilon_{\psi_1}^2(\delta^2+c \bar{d}))\frac{\alpha^{\mu \wedge k}}{\underline{p}} \bar{p}\|\tilde{x}_{0|0}\|^2,\nnum\\
&c_1=\bar{m}_1^2(\bar{c}_1^2+2\epsilon_{\psi_1}\bar{c}_1(\delta+(1+ \bar{l}\bar{c}_2) \bar{\chi})+\epsilon_{\psi_1}^2(\delta^2+c \bar{d}))\nnum\\
&\times\frac{c}{\underline{p}(1-\alpha)}+\bar{m}_1^2(\bar{c}_1^2(\delta(1+ \bar{l}\bar{c}_2) \bar{\chi}+(1+ \bar{l}\bar{c}_2)^2 \bar{\chi}^2)\nnum\\
&+\bar{c}_1c \bar{d}(\delta+(1+ \bar{l}\bar{c}_2) \bar{\chi})+\bar{r}_1 {Rank(R_{1,\mu \wedge k})})+\epsilon^2\bar{d}^2.
\end{align*}
By Markov's inequality, for $t\in[t_k,t_{k+1})$, we have
\begin{align*}
&P(\|\tilde{d}_1(\mu \wedge t)\|^2 \geq \frac{\beta_1(\|\tilde{x}_{0|0}\|^2,\mu \wedge k)+ c_1}{\gamma_3})\leq \gamma_3
\end{align*}
for any $\gamma_3 \in (0,1)$. Analogous to~\eqref{np_006}, we have the following for $t\in[t_k, t_{k+1})$
\begin{align*}
&P(\|\tilde{d}_1(t)\|^2 < \frac{\beta_1(\|\tilde{x}_{0|0}\|^2,k)+ c_1}{\gamma_3})\nnum\\
&\geq
P(\|\tilde{d}_1(\mu \wedge t)\|^2 < \frac{\beta_1(\|\tilde{x}_{0|0}\|^2,\mu \wedge k)+ c_1}{\gamma_3}| \mu > k)\nnum\\
&\quad \times P(\mu>k)
\end{align*}
\begin{align*}
&= P(\|\tilde{d}_1(\mu \wedge t)\|^2 < \frac{\beta_1(\|\tilde{x}_{0|0}\|^2,\mu \wedge k)+ c_1}{\gamma_3})\nnum\\
&\quad - P(\|\tilde{d}_1(\mu \wedge t)\|^2 < \frac{\beta_1(\|\tilde{x}_{0|0}\|^2,\mu \wedge k)+ c_1}{\gamma_3}| \mu \leq k)\nnum\\
&\quad \times P(\mu \leq k)\nnum\\
&\geq 1 - \gamma_3-\gamma_1 \geq \gamma
\end{align*}
for some $\gamma \in (0,1)$.
By applying Minkowski inequality to $\frac{\beta_1(\|\tilde{x}_{0|0}\|^2,k)+ c_1}{\gamma_3}$, we obtain PESp-like property for $\tilde{d}_1(t)$.
%where corresponding constants can be obtained as we did for $\tilde{x}_{k|k}$ in Claim B.5.
The proof of PESp-like property for $\tilde{d}_2(t)$ is similar to that for $\tilde{d}_1(t)$. We omit its details.
Thus, we complete the proof.
%Thus, we have PESp-like property for $\tilde{d}(t) = G^{\dagger}(t)(G_1(t)\tilde{d}_1(t)+G_2(t)\tilde{d}_2(t))$.
\oprocend

\subsubsection{Proof of Lemma~\ref{Bounded_2}}\label{AP_CD_Th_proof2}
To prove the statement, we first formally establish the equivalence of the NISE and the extended Kalman filter by expressing the attack vector estimates as functions of state estimates. Due to such connection, we apply existing results on the analysis of the extended Kalman filer~\cite{frogerais2012various} to the NISE to prove the rest of the part.

By expressing $\tilde{d}_{1,k-1,p}$ and $\tilde{d}_{2,k-1,p}$ as functions of $\tilde{x}_{k-1|k-1,p}$, 
the state estimation error update of the NISE can be found by~\eqref{CD004}.
The estimation update rule of~\eqref{CD004} is identical to that of the following continuous-discrete extended Kalman filtering problem
\begin{align}
\dot{x}(t) &= \bar{f}(x(t),u(t),\bar{w}'(t),t)\nnum\\
z_{3,k} &= \bar{h}(x_k,u_k,v_k',t_k)
\label{CD_equiv}
\end{align}
where its linearized system is given by
\begin{align*}
x_{k+1} &= x_k+ \epsilon \bar{f}(x_k,u_k,\bar{w}_k',t_k) +\epsilon \bar{\rho}_{k} \nnum\\
&\simeq x_k+\epsilon (\bar{A}_k x_k + \bar{B}_ku_k+ \bar{\rho}_{k})+\bar{w}_k\nnum\\
z_{3,k}&\simeq C_{3,k}x_k+D_{3,k}u_k+v_{3,k}
\end{align*}
as shown in Claim C.1; i.e., the two problems are equivalent to each other.

\textbf{Claim C.1:}
Under Assumption~\ref{CD_asm_m},
state and error covariance update rule of the continuous-discrete extended Kalman filter problem~\eqref{CD_equiv} is identical to that described in~\eqref{CD004} for the system~\eqref{CD001a}. Moreover, the optimal estimate gains are identical.

\textbf{PROOF.} 
First consider~\eqref{CD_equiv}. State prediction can be
\begin{align*}
\dot{\hat{x}}(t) = \bar{f}(\hat{x}(t),u(t),0,t)
\end{align*}
for $t \in (t_{k-1},t_k]$ to have $\hat{x}_{k|k-1}$ at $t=t_k$.
Or equivalently,
$
\hat{x}_{k|k-1} = \hat{x}_{k-1|k-1}+\epsilon \bar{f}(\hat{x}_{k-1|k-1},u_{k-1},0,t_{k-1}) +\epsilon\bar{\rho}_{k-1}.
$
Its error dynamic is given by, linearizing $f$,
\begin{align*}
\tilde{x}_{k|k-1,p} &= (I+\epsilon \bar{A}_{k-1})\tilde{x}_{k-1|k-1,p}+\bar{w}_{k-1}+ \bar{\phi}_{k-1}.
\end{align*}
State estimate is
\begin{align*}
\hat{x}_{k|k} &= \hat{x}_{k|k-1}+L_{k}(z_{3,k} - \bar{h}(x_{k|k-1},u_k,0,t_k))
\end{align*}
with its error dynamic
\begin{align*}
&\tilde{x}_{k|k,p} =(I-L_{k}C_{3,k})\tilde{x}_{k|k-1,p} -L_k v_{3,k}-L_k\psi_{3,k}\nnum\\
&=\bar{A}_{k-1}\tilde{x}_{k-1|k-1,p}+\bar{w}_{k-1}+\bar{\phi}_{k-1}+\epsilon\bar{\rho}_{k-1}\nnum\\
&-L_k (C_{3,k}(\bar{A}\tilde{x}_{k-1|k-1,p}+\bar{w}_{k-1}+\bar{\phi}_{k-1}+\epsilon\bar{\rho}_{k-1})+v_{3,k}\nnum\\
&+\psi_{3,k}).
\end{align*}
which is identical to~\eqref{CD004}. Once the state update rule~\eqref{CD004} is obtained, then error covariance update can be found by the same argument for the NISE. Moreover, it should be emphasized that we need to solve the same unconstrained optimization problem to find the optimal gains $L_k$ for the both extended Kalman filer and NISE.
\oprocend

Assumption~\ref{asm_1} implies that
$\|\bar{A}_{k} \|\leq \bar{a} $, 
$\underline{q} I \leq \bar{Q}_k$ and uniform observability of the pair $(C_{3,k},\bar{A}_k)$.
The above conditions with Assumptions~\ref{asm_1}, and~\ref{asm_2} suffice the condition of Theorem 4.5  in~\cite{kluge2010stochastic} (for the extended Kalman filter) which ensures the existence of constants $\underline{p}$ and $\bar{p}$ such that $\underline{p} I\leq P_k^x={\mathbb E}[\tilde{x}_k\tilde{x}_k^{T}] \leq \bar{p}I$ for all $k \geq 0$. We could apply the existing result because the update rule of the extended Kalman filter is equivalent to the NISE by Claim C.1.
%By Proposition 5.3 in~\cite{terrell2009stability}, $(C_{1,k-1},A_{k-1}- G_{k-1} M_{k-1} C_{1,k-1})$ is observable if $(C_{1,k-1},A_{k-1})$ is observable. 
\oprocend

The equivalence to the extended Kalman filter is achieved by the orthogonal output decompositions~\eqref{CD002} where their noise vectors are uncorrelated with each others.
This leads to the fact that process noise $\bar{w}_{k-1}$ and measurement noise $v_{3,k}$ in~\eqref{CD004} are uncorrelated as desired in the extended Kalman filtering problem.

\section{Conclusion}
We formulate the attack-resilient estimation of a class of switched nonlinear stochastic systems as the problem of joint estimation of the states, attack vectors and modes.
%We design the attack-resilient state, attack vector, and mode estimator for the estimation problem.
The proposed estimator, the NISME, consists of multiple NISE and a mode estimator. 
Each NISE is able to generate state and attack estimates for a particular mode and the mode estimator chooses the most likely one. Lastly, the NISME uses the estimates of the selected mode as outputs.
We formally analyze the stability of estimation errors in probability for the proposed estimator associated with the true mode under the time-invariant hidden mode. We propose a way to alleviate computational complexity by reducing the number of modes.
The estimator performance for time-varying modes with a regular mode set and a reduced mode set is validated by the numerical simulations on the IEEE 68-bus test system.

\end{document}